\documentclass[11pt]{article}
\usepackage[utf8]{inputenc}
\usepackage[T1]{fontenc}
\usepackage{amsmath,amssymb,amsfonts,amsthm,lmodern}
\usepackage{mathtools}
\usepackage{geometry}
\usepackage[hidelinks,
pdftitle={Extremals and normalized solutions for a critical weighted Choquard equation},
pdfauthor={},
pdfkeywords={Caffarelli--Kohn--Nirenberg equation; Stein--Weiss convolution; extremal function; global compactness; normalized solution}]{hyperref}
\theoremstyle{plain}
\newtheorem{theorem}{Theorem}[section]
\newtheorem{proposition}[theorem]{Proposition}
\newtheorem{lemma}[theorem]{Lemma}
\newtheorem{corollary}[theorem]{Corollary}

\theoremstyle{definition}

\theoremstyle{remark}
\newtheorem{remark}[theorem]{Remark}

\numberwithin{equation}{section}

\title{Extremals and Normalized Solutions for a Critical Weighted Choquard Equation}
\author{  }
\date{}

\begin{document}
\author{Wenjing Chen\footnote{Corresponding author.}\ \footnote{E-mail address:\, {\tt wjchen@swu.edu.cn} (W. Chen), {\tt zqzhamath@163.com} (Z. Zha).}\ \ and Ziqiang Zha\\
	\footnotesize  School of Mathematics and Statistics, Southwest University, Chongqing, 400715, P.R. China}
\maketitle

\begin{abstract}
	This paper is concerned with the existence and multiplicity of normalized solutions to a weighted Choquard equation involving a critical Stein--Weiss exponent
	\[
	\left\{
	\begin{aligned}
		-\operatorname{div}(|x|^{-2a}\nabla u)
		&=\lambda |x|^{-2a}u
		+\beta |x|^{-bq}|u|^{q-2}u\\
		&\quad
		+|x|^{-bp}
		\left(
		\int_{\mathbb R^N}
		\frac{|u(y)|^p}
		{|x-y|^{\mu}|y|^{bp}}\,dy
		\right)
		|u|^{p-2}u
		&&\text{in }\mathbb R^N,\\
		\int_{\mathbb R^N}|x|^{-2a}u^2\,dx&=\rho^2,
	\end{aligned}
	\right.
	\]
	where $N\geq3$, $a<(N-2)/2$, $a<b<a+1$,
	$2<q<2^\sharp$, $0<\mu<N$, $p:=(2N-\mu)/(N-2d)$ is the upper Stein-Weiss
	critical exponent, $d=1-b+a$ and $2^\sharp:={2N}/({N-2-2a+2b})$ is the Caffarelli-Kohn-Nirenberg critical exponent. The weighted mass $\rho>0$ is prescribed,
	$\beta>0$, and $\lambda\in\mathbb R$
	is an unknown Lagrange multiplier.
	
	Assuming $p>2$, we establish the existence and multiplicity of normalized solutions across different mass regimes. In the mass-subcritical regime, we establish the existence of a negative-energy normalized ground state for sufficiently small $\beta$. Under suitable assumptions, a second positive normalized solution at a positive energy level is also obtained. The nonlocal interaction provides an alternative sufficient condition for the second solution, covering cases in which the local perturbation estimate does not apply. Furthermore, we prove the existence of normalized ground states in both the mass-critical and mass-supercritical regimes under explicit parameter conditions. The proof relies on a global compactness description for constrained Palais--Smale sequences, together with the existence and sharp decay asymptotics of positive extremals for the limiting critical problem. By designing minimax test functions via suitably truncated extremals, we dispense with the restrictive assumption of finite weighted mass for the limiting profiles.
	\medskip
	
\end{abstract}


\section{Introduction and main results}\label{sec:introduction}

Standing waves with prescribed mass have been extensively studied
for nonlinear Schr\"odinger equations. In this setting, the mass
is fixed, while the frequency is an unknown Lagrange multiplier.
Cazenave and Lions \cite{CazenaveLions1982} used constrained
minimization to obtain standing waves and prove the orbital stability of these waves. Jeanjean \cite{Jeanjean1997} developed a minimax
approach to solutions with prescribed norms. For equations with
combined power nonlinearities, Soave
\cite{SoaveJDE2020,SoaveJFA2020} studied normalized ground states
in the Sobolev subcritical and critical cases. Further existence
and multiplicity results in the Sobolev critical case were
obtained by Jeanjean and Le \cite{JeanjeanLe2022} and Wei and Wu
\cite{WeiWu2022NLS}.

Equations involving singular potentials or weighted nonlinearities
have also received considerable attention. For local critical
equations with an inverse-square potential, Terracini
\cite{Terracini1996} studied the existence and qualitative
properties of solutions. In the prescribed-mass setting, Cardoso,
de Oliveira, and Miyagaki \cite{CardosoOliveiraMiyagaki2025}
considered equations with Hardy--Sobolev nonlinearities.
Besides equations with weights in the lower-order terms,
one can consider equations with a weighted principal part,
such as $-\operatorname{div}(|x|^{-2a}\nabla u)$.
The Caffarelli--Kohn--Nirenberg inequality \cite{CKN1984}
provides a basic estimate for such problems. More precisely,
assume that
\begin{equation}\label{eq:assumption-ab}
	N\geq3,\qquad
	a<\frac{N-2}{2},\qquad
	a<b<a+1.
\end{equation}
Then there exists a best constant $S(a,b)>0$ such that
\begin{equation}\label{eq:critical-CKN-intro}
	S(a,b)
	\left(
	\int_{\mathbb R^N}
	|x|^{-b2^\sharp}|u|^{2^\sharp}\,dx
	\right)^{2/2^\sharp}
	\leq
	\int_{\mathbb R^N}|x|^{-2a}|\nabla u|^2\,dx,
	\qquad u\in C_c^\infty(\mathbb R^N),
\end{equation}
where
\[
2^\sharp=\frac{2N}{N-2-2a+2b}.
\]
The best constants and the existence and symmetry properties
of extremals were studied by Catrina and Wang
\cite{CatrinaWang2001} and Felli and Schneider
\cite{FelliSchneider2003}. The cylindrical formulation and
rigidity results were developed in
\cite{DolbeaultEstebanLoss2016}. For quantitative stability
and the existence of minimizers for the optimal stability
constant, we refer to
\cite{WeiWu2022CKN,FrankPeteranderl2024,DengTianWei2025}.
In this weighted setting, the prescribed-mass problem was
studied by Goel and Rai \cite{GoelRai2026}. They considered
\begin{equation}\label{eq:GoelRai}
	\left\{
	\begin{aligned}
		-\operatorname{div}(|x|^{-2a}\nabla u)
		&=
		\lambda |x|^{-2a}u
		+\beta |x|^{-bq}|u|^{q-2}u
		+|x|^{-b2^\sharp}|u|^{2^\sharp-2}u
		&&\text{in }\mathbb R^N,\\
		\int_{\mathbb R^N}|x|^{-2a}u^2\,dx
		&=\rho^2,
	\end{aligned}
	\right.
\end{equation}
where $0<a<(N-2)/2$, $a<b<a+1$, $2<q<2^\sharp$,
$\beta,\rho>0$, and $\lambda\in\mathbb R$ is an unknown
Lagrange multiplier. Within a radial energy space, they
obtained two positive normalized solutions in the
mass-subcritical regime and normalized ground states in
the mass-critical and mass-supercritical regimes,
under suitable conditions on the parameters.

Normalized solutions have also been studied for equations of critical Choquard type. The nonlocal energy is estimated by the
Hardy--Littlewood--Sobolev inequality \cite{Lieb1983}.
The existence and qualitative properties of Choquard ground
states were studied by Moroz and Van Schaftingen
\cite{MorozVanSchaftingen2013,MorozVanSchaftingen2017}.
For equations with an upper critical Choquard term and a local
power perturbation, Li \cite{Li2022Standing,Li2023Normalized}
established multiplicity results in the mass-subcritical case
and existence and nonexistence results for normalized ground
states in the mass-critical and mass-supercritical cases.

When singular weights are included in the Choquard term,
the corresponding energy estimates rely on the Stein--Weiss
inequality \cite{SteinWeiss1958}. Without a prescribed-mass
constraint, Du, Gao, and Yang \cite{DuGaoYang2022} proved
the attainment of the critical Stein--Weiss quotient, while
Li, Yang, and Zhou \cite{LiYangZhou2022} studied qualitative
properties of solutions and established Liouville-type results for Stein--Weiss equations. In dimension two,
Alves and Shen \cite{AlvesShen2023} studied Stein--Weiss
problems with critical exponential growth. For Choquard
equations with a Hardy potential, Guo, Tang, and Zhang
\cite{GuoTangZhang2022} obtained ground states in the presence
of a local Sobolev critical term. Guo and Tang
\cite{GuoTang2023} established existence and sharp two-sided
asymptotic estimates for positive solutions of the critical
problem. In the normalized
setting, Guo and Zhang \cite{GuoZhang2025} studied Choquard
equations with Hardy terms and local perturbations under
different assumptions on the potential.

These developments suggest considering a nonlocal counterpart of
\eqref{eq:GoelRai}, in which the local critical power is replaced
by an upper critical Stein--Weiss term. We therefore study the
prescribed-mass problem
\begin{equation}\label{eq:problem}
	\left\{
	\begin{aligned}
		-\operatorname{div}(|x|^{-2a}\nabla u)
		&=
		\lambda |x|^{-2a}u
		+\beta |x|^{-bq}|u|^{q-2}u\\
		&\quad
		+|x|^{-bp}
		\left(
		\int_{\mathbb R^N}
		\frac{|u(y)|^p}
		{|x-y|^\mu |y|^{bp}}\,dy
		\right)
		|u|^{p-2}u
		&&\text{in }\mathbb R^N,\\
		\int_{\mathbb R^N}|x|^{-2a}u^2\,dx
		&=\rho^2,
	\end{aligned}
	\right.
\end{equation}
where $\rho>0$, $\beta>0$, and $\lambda\in\mathbb R$
is an unknown Lagrange multiplier.

Under assumption \eqref{eq:assumption-ab}, we set
\[
\kappa=N-2-2a>0,\qquad
\alpha=b-a\in(0,1),\qquad
d=1-\alpha=1+a-b\in(0,1).
\]
The critical exponents are
\begin{equation}\label{eq:critical-exponents-intro}
	2^\sharp=\frac{2N}{N-2d},\qquad
	0<\mu<N,\qquad
	p=\frac{2N-\mu}{N-2d}>1.
\end{equation}
For the local perturbation, we assume $2<q<2^\sharp$.
Notice that $2^\sharp<2^*=2N/(N-2)$. The following decay and normalized existence results are proved under the additional assumption
\begin{equation}\label{eq:decay-range}
	0<\mu<\min\{N,4d\},
\end{equation}
which is equivalent to $p>2$. Unlike the framework in \cite{GoelRai2026}, we investigate a critical nonlocal term, relax the parameter restriction to include $a\leq0$, and establish our results in the full energy space.

Problem \eqref{eq:problem} can be rewritten as a
Schr\"odinger equation with a Hardy potential through the
change of variables $v(x)=|x|^{-a}u(x)$. Integration by
parts gives
\[
\int_{\mathbb R^N}|x|^{-2a}|\nabla u|^2\,dx
=
\int_{\mathbb R^N}|\nabla v|^2\,dx
-h_a\int_{\mathbb R^N}\frac{|v|^2}{|x|^2}\,dx,
\qquad
h_a=a(N-2-a)<\frac{(N-2)^2}{4}.
\]
This identity holds first for smooth functions supported
away from the origin and extends to the energy space by
density and Hardy's inequality. Thus, \eqref{eq:problem}
is equivalent to
\begin{equation}\label{eq:hardy-formulation}
	\left\{
	\begin{aligned}
		-\Delta v-\frac{h_a}{|x|^2}v
		&=
		\lambda v
		+\beta|x|^{-\alpha q}|v|^{q-2}v\\
		&\quad
		+|x|^{-\alpha p}
		\left(
		\int_{\mathbb R^N}
		\frac{|v(y)|^p}
		{|x-y|^\mu|y|^{\alpha p}}\,dy
		\right)
		|v|^{p-2}v
		&&\text{in }\mathbb R^N,\\
		\int_{\mathbb R^N}|v|^2\,dx
		&=\rho^2.
	\end{aligned}
	\right.
\end{equation}
The nonlocal weight $\zeta=\alpha p>0$ satisfies
\[
(N-2)p+2\zeta+\mu=2N,
\]
which is the critical scaling relation for the Stein--Weiss
term. When $a<0$, one has $h_a<0$, so the inverse-square
term makes a positive contribution to the energy.
Schwarz rearrangement need not decrease the full quadratic
form in this case and therefore does not justify restricting
the variational problem to radial functions.

The critical nonlocal term leads to a loss of compactness
in the weighted energy space. The concentration--compactness
principle of Lions \cite{Lions1984} and the global compactness
theorem of Struwe \cite{Struwe1984} provide the classical
framework, while corresponding nonlocal compactness results
for critical Choquard equations were established in
\cite{GaoDaSilvaYangZhou2020,HuangLongXiaZheng2019}.
In the present work, we establish profile decompositions
for both the limiting critical functional and constrained
Palais--Smale sequences in the full weighted energy space.
We also prove that the critical best constant is attained
and determine the precise decay of positive extremals.
These results are used to obtain normalized ground states
and a second positive solution under conditions arising
from estimates of the local perturbation or the nonlocal
interaction. The comparison arguments use truncated
extremals and do not require $\kappa>2$, the condition
under which the critical extremals have finite mass.

For convenience, set
\[
\begin{aligned}
	A(u)&=\int_{\mathbb R^N}|x|^{-2a}|\nabla u|^2\,dx,\qquad
	M(u)=\int_{\mathbb R^N}|x|^{-2a}|u|^2\,dx,\\
	B_q(u)&=\int_{\mathbb R^N}|x|^{-bq}|u|^q\,dx,\qquad
	\mathcal D(u)=\int_{\mathbb R^N}\int_{\mathbb R^N}
	\frac{|u(x)|^p|u(y)|^p}{|x|^{bp}|x-y|^{\mu}|y|^{bp}}\,dx\,dy.
\end{aligned}
\]
For $1\leq s<\infty$, we write $L^s(\mathbb R^N;|x|^{-bs})$ for the weighted Lebesgue space with norm
\[
\|u\|_{L^s(\mathbb R^N;|x|^{-bs})}
=\left(\int_{\mathbb R^N}|x|^{-bs}|u|^s\,dx\right)^{1/s}.
\]
Let $D_a^{1,2}(\mathbb R^N)$ be the completion of $C_c^\infty(\mathbb R^N)$ under the norm $\|u\|_a=A(u)^{1/2}$, and set
\[
X=D_a^{1,2}(\mathbb R^N)\cap L^2(\mathbb R^N;|x|^{-2a}),
\qquad \|u\|_X^2=A(u)+M(u).
\]
An application of the classical Hardy--Littlewood--Sobolev inequality \cite{Lieb1983} and the Caffarelli--Kohn--Nirenberg inequality \eqref{eq:critical-CKN-intro} ensures that $\mathcal D(u)$ is well-defined on $D_a^{1,2}(\mathbb R^N)$ and satisfies the bound $\mathcal D(u)^{1/p}\leq C A(u)$.

The energy functional and the mass constraint are
\[
E_{\beta}(u)=\frac12A(u)-\frac{\beta}{q}B_q(u)-\frac{1}{2p}\mathcal D(u),\qquad S_{\rho}=\{u\in X:M(u)=\rho^2\}.
\]
A normalized solution is a critical point of $E_\beta|_{S_\rho}$. We call it a normalized ground state if it has the least energy among all critical points of $E_\beta|_{S_\rho}$.

The critical best constant and the associated energy threshold are
\[
S_{\mu}(a,b)=\inf_{u\in D_a^{1,2}\setminus\{0\}}\frac{A(u)}{\mathcal D(u)^{1/p}}>0,
\qquad
T_{\mu,a,b}=\frac{p-1}{2p}S_{\mu}(a,b)^{p/(p-1)}.
\]

The mass constraint is invariant under the scaling
\[
(t\star u)(x)=e^{(N-2a)t/2}u(e^t x),\qquad t\in\mathbb R.
\]
A direct computation shows that
\[
A(t\star u)=e^{2t}A(u),\qquad
B_q(t\star u)=e^{r_qt}B_q(u),\qquad
\mathcal D(t\star u)=e^{2pt}\mathcal D(u),
\]
where
\[
\delta_q=\frac{(N-2a+2b)q-2N}{2q},
\qquad r_q=q\delta_q,
\qquad \gamma_q=q-r_q=N-\frac{q(N-2d)}2\in(0,2).
\]
The associated fiber map is

\[\Phi_{\beta,u}(t)=E_\beta(t\star u)
=\frac12e^{2t}A(u)-\frac{\beta}{q}e^{r_qt}B_q(u)
-\frac1{2p}e^{2pt}\mathcal D(u),
\]
Then
\[
\begin{aligned}
	\Phi_{\beta,u}'(t)
	&=e^{2t}A(u)-\beta\delta_qe^{r_qt}B_q(u)-e^{2pt}\mathcal D(u),\\
	\Phi_{\beta,u}''(t)
	&=2e^{2t}A(u)-\beta q\delta_q^2e^{r_qt}B_q(u)-2pe^{2pt}\mathcal D(u).\\
	P_\beta(u)
	&=\Phi_{\beta,u}'(0)
	=A(u)-\beta\delta_qB_q(u)-\mathcal D(u).
\end{aligned}
\]
The Pohozaev manifold is
\[
\mathcal M_{\rho,\beta}=\{u\in S_{\rho}:P_{\beta}(u)=0\}.
\]
For $u\in\mathcal M_{\rho,\beta}$, set
\begin{equation}\label{eq:Q-definition}
	Q_{\beta}(u)=\Phi_{\beta,u}''(0)
	=(2-r_q)\beta\delta_qB_q(u)-(2p-2)\mathcal D(u).
\end{equation}
When $r_q<2$, we decompose $\mathcal M_{\rho,\beta}$ into
\[
\begin{aligned}
	\mathcal M_{\rho,\beta}^{+}
	&=\{u\in\mathcal M_{\rho,\beta}:Q_\beta(u)>0\},\\
	\mathcal M_{\rho,\beta}^{-}
	&=\{u\in\mathcal M_{\rho,\beta}:Q_\beta(u)<0\},\\
	\mathcal M_{\rho,\beta}^{0}
	&=\{u\in\mathcal M_{\rho,\beta}:Q_\beta(u)=0\}.
\end{aligned}
\]

The mass-critical exponent\[q_c=\frac{2N+4}{N-2a+2b} = \frac{2(N+2)}{N+2-2d}\]is determined by the threshold condition $r_{q_c}=2$. Under the mass-preserving scaling, the relative scaling order of the local perturbation with respect to the quadratic term naturally classifies the problem into three distinct regimes: the mass-subcritical ($q < q_c$), mass-critical ($q = q_c$), and mass-supercritical ($q > q_c$) cases. However, the critical nonlocal term dominates the quadratic term
under large dilations in all three regimes, so $E_\beta$ is
unbounded below on $S_\rho$. This rules out global minimization on the mass constraint and requires a more delicate variational
approach to overcome the resulting loss of compactness.

Our first result concerns the critical extremals.

\begin{theorem}[Critical extremals and decay]\label{thm:critical-extremal}
	Assume \eqref{eq:assumption-ab} and $0<\mu<N$. The constant $S_{\mu}(a,b)$ is attained by a nonnegative function $U\in D_a^{1,2}(\mathbb R^N)$. After multiplication by a positive constant, $U$ satisfies
	\begin{equation}\label{eq:pure-critical-equation-intro}
		-\operatorname{div}(|x|^{-2a}\nabla U)
		=|x|^{-bp}\left(\int_{\mathbb R^N}
		\frac{|U(y)|^p}{|x-y|^{\mu}|y|^{bp}}\,dy\right)|U|^{p-2}U
	\end{equation}
	and
	\begin{equation}\label{eq:critical-normalization-intro}
		A(U)=\mathcal D(U)=S_{\mu}(a,b)^{p/(p-1)}.
	\end{equation}
	Moreover, $U>0$ almost everywhere. If $p>2$, its continuous representative is strictly positive on $\mathbb R^N\setminus\{0\}$, and there exist positive constants $c,C$ such that
	\begin{equation}\label{eq:extremal-decay-intro}
		\frac{c}{1+|x|^{\kappa}}
		\leq U(x)\leq
		\frac{C}{1+|x|^{\kappa}},
		\qquad x\in\mathbb R^N\setminus\{0\}.
	\end{equation}
	Under the same assumption $p>2$, there exist positive constants $C_+,C_-$ such that, uniformly in $\omega\in\mathbb S^{N-1}$,
	\begin{equation}\label{eq:extremal-endpoint-limits-intro}
		U(r\omega)\longrightarrow C_+\quad(r\downarrow0),
		\qquad r^\kappa U(r\omega)\longrightarrow C_-\quad(r\to\infty).
	\end{equation}
	The explicit integral representations of $C_\pm$ are given in Proposition \ref{prop:axial-decay}.
\end{theorem}

\begin{theorem}[Global compactness for the unconstrained and constrained problems]\label{thm:global-compactness-intro}
	Assume \eqref{eq:assumption-ab} and \eqref{eq:critical-exponents-intro}, and set
	\[
	J(u)=\frac12A(u)-\frac{1}{2p}\mathcal D(u),\qquad u\in D_a^{1,2}(\mathbb R^N).
	\]
	\begin{itemize}
		\item[(i)] Let $\{u_n\}$ be bounded in $D_a^{1,2}$ and satisfy
		\[
		J(u_n)\to c,\qquad J'(u_n)\to0
		\quad\text{in }(D_a^{1,2})^*.
		\]
		Then, after passing to a subsequence, there are a solution $u^0$ of \eqref{eq:pure-critical-equation-intro}, an integer $k\geq0$, nontrivial solutions $U^1,\ldots,U^k$ of the same equation, and scales $\varepsilon_n^j>0$ such that
		\begin{equation}\label{eq:profile-decomposition-intro}
			u_n=u^0+\sum_{j=1}^k(\varepsilon_n^j)^{-\kappa/2}
			U^j\left(\frac{x}{\varepsilon_n^j}\right)+o(1)
			\quad\text{in }D_a^{1,2}(\mathbb R^N),
		\end{equation}
		where
		\begin{equation}\label{eq:scale-orthogonality-intro}
			|\log\varepsilon_n^j|\to\infty,\qquad
			\left|\log\frac{\varepsilon_n^i}{\varepsilon_n^j}\right|\to\infty
			\quad(i\neq j).
		\end{equation}
		The norm and energy split:
		\begin{equation}\label{eq:profile-splitting-intro}
			A(u_n)=A(u^0)+\sum_{j=1}^kA(U^j)+o(1),
			\qquad
			c=J(u^0)+\sum_{j=1}^kJ(U^j).
		\end{equation}
		The weak limit and the profiles may change sign. Every nontrivial profile has energy at least $T_{\mu,a,b}$.
		
		\item[(ii)] Assume also $p>2$, $2<q<2^\sharp$, and $\beta,\rho>0$. Let $\{u_n\}\subset S_\rho$ be bounded in $X$ and satisfy
		\[
		E_\beta(u_n)\to c,\qquad
		\|d(E_\beta|_{S_\rho})(u_n)\|\to0,\qquad
		P_\beta(u_n)\to0.
		\]
		Then, up to a subsequence, \eqref{eq:profile-decomposition-intro} holds for some integer $k\geq0$, with nontrivial, possibly sign-changing, solutions of the limiting critical equation $U^j$ and scales $\varepsilon_n^j\to0$ satisfying \eqref{eq:scale-orthogonality-intro}. The weak limit $u^0\in X$ solves the Euler equation in \eqref{eq:problem} with a multiplier $\lambda\leq0$, and
		\[
		c=E_\beta(u^0)+\sum_{j=1}^kJ(U^j).
		\]
		If $u^0\not\equiv0$, then $M(u^0)=\rho^2$ and $\lambda<0$; if also $k=0$, then $u_n\to u^0$ in $X$. If $u^0\equiv0$, then $\lambda=0$ and either $c=0$ (when $k=0$, meaning $u_n\to0$ in $D_a^{1,2}(\mathbb R^N)$) or $c\geq T_{\mu,a,b}$ (when $k\geq1$). Moreover,
\[
\begin{aligned}
A(u_n)&=A(u^0)+\sum_{j=1}^k A(U^j)+o(1),\\
\mathcal D(u_n)&=\mathcal D(u^0)+\sum_{j=1}^k\mathcal D(U^j)+o(1),\\
B_q(u_n)&=B_q(u^0)+o(1).
\end{aligned}
\]
	\end{itemize}
\end{theorem}

Some of the normalized existence results require the additional condition
\begin{equation}\label{eq:absolute-truncation-condition}
	\kappa>\frac{2\gamma_q}{2+\gamma_q}.
\end{equation}
Let
\begin{equation}\label{C_{q_c}}
\mathcal C_{q_c}=
\sup_{u\in X\setminus\{0\}}
\frac{B_{q_c}(u)}{A(u)M(u)^{(q_c-2)/2}},
\qquad
\beta_* =\frac{q_c}{2\mathcal C_{q_c}\rho^{q_c-2}}.
\end{equation}

\begin{theorem}[Normalized solutions]\label{thm:main}
	Assume \eqref{eq:assumption-ab}, \eqref{eq:critical-exponents-intro}, \eqref{eq:decay-range}, $2<q<2^\sharp$, and $\rho>0$.
	
	\begin{itemize}
		\item[(i)] If $2<q<q_c$, there exists $\beta_0=\beta_0(N,a,b,\mu,q,\rho)>0$ such that, for $0<\beta<\beta_0$, problem \eqref{eq:problem} has a positive normalized ground state $u_{\beta}^{+}\in\mathcal M_{\rho,\beta}^{+}$ satisfying
		\[
		E_{\beta}(u_{\beta}^{+})<0,
		\qquad \lambda_{\beta}^{+}<0.
		\]
		It is a local minimizer of $E_{\beta}$ on $S_{\rho}$.
		
		\item[(ii)] If
		\begin{equation}\label{eq:main-multiplicity-range}
			2<q<q_c,\qquad
			\left(\gamma_q<\frac{\kappa}{2}\quad\text{or}\quad\mu<\kappa\right),
		\end{equation}
		then, after decreasing $\beta_0$ if necessary, for every $0<\beta<\beta_0$ problem \eqref{eq:problem} has a second positive normalized solution $u_{\beta}^{-}\in\mathcal M_{\rho,\beta}^{-}$ with
		\[
		0<E_{\beta}(u_{\beta}^{-}),
		\qquad \lambda_{\beta}^{-}<0.
		\]
		Its energy satisfies
		\[
		E_{\beta}(u_{\beta}^{-})
		<E_{\beta}(u_{\beta}^{+})+T_{\mu,a,b}.
		\]
		
		\item[(iii)] If $q=q_c$, $0<\beta<\beta_*$, and $\kappa>4d/(N+2)$, problem \eqref{eq:problem} has a positive normalized ground state $u_{\beta}$ such that
		\begin{equation}\label{eq:main-critical}
			0<E_{\beta}(u_{\beta})<T_{\mu,a,b},
			\qquad \lambda_{\beta}<0.
		\end{equation}
		
		\item[(iv)] If $q_c<q<2^{\sharp}$ and \eqref{eq:absolute-truncation-condition} holds, then for every $\beta>0$ problem \eqref{eq:problem} has a positive normalized ground state satisfying \eqref{eq:main-critical}.
	\end{itemize}
\end{theorem}

\begin{remark}\label{rem:main-thm-comments}
	The two conditions in Theorem \ref{thm:main}\textit{(ii)} stem from estimates of the local perturbation and the nonlocal interaction, respectively. The first is equivalent to
	\[
	\max\left\{2,\frac{2N-\kappa}{N-2d}\right\}<q<q_c,
	\]
	which defines a nonempty interval precisely when $\kappa>8d/(N+2-2d)$. In contrast, provided \eqref{eq:decay-range} holds, the second condition $\mu<\kappa$ applies to the full range $2<q<q_c$. This establishes existence in cases not covered by the first condition. For example, taking $N=3$, $a=0$, $b=1/2$, $\mu=1/2$, and $q=11/5$ yields $\mu<\kappa=1$, whereas $\gamma_q=4/5>\kappa/2$. Furthermore, the conditions in parts \textit{(iii)} and \textit{(iv)} are derived directly from \eqref{eq:absolute-truncation-condition}. At $q=q_c$, this reduces to $\kappa>4d/(N+2)$. Since $\gamma_q$ is strictly decreasing with respect to $q$, the condition in part \textit{(iv)} holds for every $q_c<q<2^\sharp$ whenever $\kappa\geq4d/(N+2)$. Finally, it is worth noting that Theorem \ref{thm:main} does not impose the finite mass requirement on the critical extremals. Although the decay estimate in Theorem \ref{thm:critical-extremal} implies that a critical extremal belongs to $X$ if and only if $\kappa>2$ under \eqref{eq:decay-range}, we overcome this restriction by constructing our test functions via truncated critical extremals.
\end{remark}

The paper is organized as follows. Section \ref{sec:framework} collects preliminary estimates and develops the variational framework. In Section \ref{sec:critical-compactness}, we establish the existence and asymptotic decay of critical extremals and prove the global compactness theorems. Section \ref{sec:normalized-solutions} is devoted to the normalized problem, where we establish the existence of ground states and obtain multiplicity results under the stated assumptions.

\section{Preliminaries and variational setting}\label{sec:framework}
We begin with basic weighted estimates, followed by the
Emden--Fowler transformation and the nonlocal
Br\'{e}zis--Lieb splitting. We then derive the Pohozaev
identities on the cylinder and use mass-preserving fiber
maps to characterize the geometry of the Pohozaev manifold.

\subsection{Basic weighted estimates}\label{sec:weighted-estimates}

For $2\leq s\leq2^{\sharp}$, put
\[
\delta_s=\frac{(N-2a+2b)s-2N}{2s}.
\]

\begin{lemma}[Weighted interpolation]\label{lem:weighted-interpolation}\cite{CKN1984}
	For every $s\in(2,2^{\sharp}]$, there exists $C_{a,b,s}>0$ such that
	\begin{equation}\label{eq:weighted-interpolation}
		\int_{\mathbb R^N}|x|^{-bs}|u|^s\,dx
		\leq C_{a,b,s}A(u)^{s\delta_s/2}M(u)^{s(1-\delta_s)/2},
		\qquad u\in X.
	\end{equation}
	At $s=2^{\sharp}$,
	\begin{equation}\label{eq:critical-embedding}
		\left(\int_{\mathbb R^N}|x|^{-b2^{\sharp}}|u|^{2^{\sharp}}\,dx\right)^{2/2^{\sharp}}
		\leq S(a,b)^{-1}A(u).
	\end{equation}
\end{lemma}

\begin{proposition}[Compact weighted embedding]\label{prop:compact-embedding}
	For every $s\in(2,2^{\sharp})$, the embedding
	\[
	X\hookrightarrow L^s(\mathbb R^N;|x|^{-bs})
	\]
	is compact.
\end{proposition}
\begin{proof}
Let $u_n \rightharpoonup u$ in $X$ and set $v_n := u_n - u$, which yields $A(v_n) + M(v_n) \leq C$. By Rellich's compactness theorem, we have $v_n \to 0$ strongly in $L^s(\Omega_{r,R}; \vert{}x\vert{}^{-bs})$ on any annulus $\Omega_{r,R} = \{x \in \mathbb{R}^N : r < \vert{}x\vert{} < R\}$. Fixing $\vartheta \in (0, 1)$ with $1/s = \vartheta/2 + (1-\vartheta)/2^\sharp$, by H\"older's inequality and \eqref{eq:critical-embedding}, we obtain
\[
\begin{aligned}
\int_{|x|<r}|x|^{-bs}|v_n|^s\,dx
&\leq |B_r|^{1-s/2^\sharp}
\bigl\||x|^{-b}v_n\bigr\|_{2^\sharp}^{s}
\leq Cr^{N(1-s/2^\sharp)},\\
\int_{|x|>R}|x|^{-bs}|v_n|^s\,dx
&\leq R^{-(b-a)\vartheta s}M(v_n)^{\vartheta s/2}
\bigl\||x|^{-b}v_n\bigr\|_{2^\sharp}^{(1-\vartheta)s}
\leq CR^{-(b-a)\vartheta s}.
\end{aligned}
\]
The second estimate follows from the fact that
$\int_{|x|>R}|x|^{-2b}|v_n|^2\,dx\leq R^{-2(b-a)}M(v_n)$.
Combining these estimates and letting $r\downarrow0$ and $R\to\infty$ proves the assertion.
\end{proof}

Let
\[
r_{\mu}=\frac{2N}{2N-\mu}.
\]
Then $pr_{\mu}=2^{\sharp}$. Applying the Hardy--Littlewood--Sobolev inequality to $|x|^{-bp}|u|^p\in L^{r_\mu}(\mathbb R^N)$ and \eqref{eq:critical-embedding} gives
\begin{equation}\label{eq:D-continuity}
	\mathcal D(u)\leq C
	\left(\int_{\mathbb R^N}|x|^{-b2^{\sharp}}|u|^{2^{\sharp}}\,dx\right)^{2p/2^{\sharp}}
	\leq C A(u)^p.
\end{equation}

\subsection{Cylindrical formulation}\label{sec:cylinder}
To handle the singular weights and the dilation invariance of the critical problem, we introduce the Emden--Fowler transformation. It converts the weighted equation into an elliptic problem on the cylinder $\mathcal C=\mathbb R\times\mathbb S^{N-1}$, where critical dilations correspond to axial translations. This formulation is convenient for analyzing the asymptotic decay of critical extremals and the loss of compactness.

For $x=e^{-t}\omega$ with $(t,\omega)\in\mathcal C=\mathbb R\times\mathbb S^{N-1}$, define
\[
w(t,\omega)=e^{-\kappa t/2}u(e^{-t}\omega),
\qquad\text{or equivalently,}\qquad
u(x)=|x|^{-\kappa/2}w(-\log|x|,x/|x|).
\]
Setting
\[
L_\kappa=-\partial_{tt}-\Delta_\omega+\frac{\kappa^2}{4},
\]
A direct calculation gives
\begin{equation}\label{eq:cylinder-A}
	A(u)=\int_{\mathcal C}
	\left(|\partial_tw|^2+|\nabla_\omega w|^2
	+\frac{\kappa^2}{4}w^2\right)\,dt\,d\omega.
\end{equation}
We equip $H^1(\mathcal C)$ with the equivalent norm
\[
\|w\|_\kappa^2
:=\int_{\mathcal C}
\left(|\partial_tw|^2+|\nabla_\omega w|^2
+\frac{\kappa^2}{4}w^2\right)\,dt\,d\omega.
\]
Since $\kappa>0$, this norm is equivalent to the
standard $H^1(\mathcal C)$ norm. Thus the
Emden--Fowler transformation is an isomorphism between
$D_a^{1,2}(\mathbb R^N)$ and $H^1(\mathcal C)$.
Under this transformation, the mass and the local term become
\begin{equation}\label{eq:cylinder-mass-local}
	M(u)=\int_{\mathcal C}e^{-2t}w^2\,dt\,d\omega,
	\qquad
	B_q(u)=\int_{\mathcal C}e^{-\gamma_qt}|w|^q\,dt\,d\omega.
\end{equation}

Set
\[
\mathcal K_{\mu}(s,\omega,\omega')
=\left[2(\cosh s-\omega\cdot\omega')\right]^{-\mu/2}.
\]
The critical relation $p(2b+\kappa)+\mu=2N$ shows that
\[
\mathcal D(u)=\int_{\mathcal C}\int_{\mathcal C}
\mathcal K_{\mu}(t-s,\omega,\omega')
|w(t,\omega)|^p|w(s,\omega')|^p
\,dt\,d\omega\,ds\,d\omega'.
\]
For a measurable function $f$ for which the following integral is defined, write
\[
(\mathcal K_\mu*f)(t,\omega)
:=\int_{\mathcal C}
\mathcal K_\mu(t-s,\omega,\omega')f(s,\omega')\,ds\,d\omega'.
\]
Here $*$ denotes this integral on the cylinder, and $d\omega$ and $d\omega'$ denote surface measure on $\mathbb S^{N-1}$. We also set
\[
\mathfrak B_{\mathcal C}(f,g)
:=\int_{\mathcal C}(\mathcal K_\mu*f)(t,\omega)g(t,\omega)\,dt\,d\omega,
\qquad
\mathcal D_{\mathcal C}(w):=\mathfrak B_{\mathcal C}(|w|^p,|w|^p).
\]
Thus $\mathcal D_{\mathcal C}(w)=\mathcal D(u)$. In particular,
\[
\mathcal D_{\mathcal C}'(w)[\varphi]
=2p\int_{\mathcal C}(\mathcal K_\mu*|w|^p)|w|^{p-2}w\varphi\,dt\,d\omega.
\]
The kernel $\mathcal K_\mu$ is invariant under simultaneous axial translations. The Euclidean critical dilation
\[
u_{\varepsilon}(x)=\varepsilon^{-\kappa/2}u(x/\varepsilon)
\]
corresponds to the axial translation $w_\varepsilon(t,\omega)=w(t+\log\varepsilon,\omega)$.

For $f\in L^{r_\mu}(\mathcal C)$, set
$\widetilde f(e^{-t}\omega)=e^{(N-\mu/2)t}f(t,\omega)$ and define $\widetilde g$ analogously.
Since $(N-\mu/2)r_\mu=N$, a change of variables gives
\[
\|\widetilde f\|_{L^{r_\mu}(\mathbb R^N)}=\|f\|_{L^{r_\mu}(\mathcal C)},\qquad
\mathfrak B_{\mathcal C}(f,g)
=\iint_{\mathbb R^N\times\mathbb R^N}
\frac{\widetilde f(x)\widetilde g(y)}{|x-y|^\mu}\,dx\,dy.
\]
The Euclidean HLS inequality therefore yields
$|\mathfrak B_{\mathcal C}(f,g)|\leq C\|f\|_{L^{r_\mu}(\mathcal C)}\|g\|_{L^{r_\mu}(\mathcal C)}$.
Taking $f=g=|w|^p$ and noting that $2^{\sharp}<2^*$, the Sobolev embedding on $\mathcal C$ implies
\begin{equation}\label{eq:D-cylinder-bound}
	\mathcal D_{\mathcal C}(w)
	\leq C\|w\|_{L^{2^{\sharp}}(\mathcal C)}^{2p}
	\leq C\|w\|_{H^1(\mathcal C)}^{2p}.
\end{equation}

\begin{lemma}[Nonlocal Br\'{e}zis--Lieb splitting]\label{lem:nonlocal-BL}
	Let $u_n\rightharpoonup u$ in $D_a^{1,2}(\mathbb R^N)$ and assume $u_n\to u$ almost everywhere in $\mathbb R^N$. Put $v_n=u_n-u$. Then
	\begin{equation}\label{eq:D-BL-value}
		\mathcal D(u_n)=\mathcal D(u)+\mathcal D(v_n)+o(1).
	\end{equation}
	Moreover,
	\begin{equation}\label{eq:D-BL-derivative}
		\mathcal D'(u_n)-\mathcal D'(u)-\mathcal D'(v_n)\to0
		\quad\text{in }\bigl(D_a^{1,2}(\mathbb R^N)\bigr)^*.
	\end{equation}
	Finally, for any fixed $\varphi\in D_a^{1,2}(\mathbb R^N)$, we have
	\[
	\mathcal D'(v_n)[\varphi]\to0.
	\]
\end{lemma}

\begin{proof}
	Let $w_n$ and $w$ be the Emden--Fowler transforms of $u_n$ and $u$, respectively. Write $w_n=w+z_n$ and set
	\[
	\ell=2^\sharp,\qquad r=r_\mu=\frac\ell p,
	\qquad r'=\frac{2N}\mu,\qquad h=\frac\ell{p-1},
	\qquad \ell'=\frac\ell{\ell-1}.
	\]
	Thus $1/\ell'=1/r'+1/h$. By the cylindrical HLS estimate and duality, we obtain
	\begin{equation}\label{eq:BL-operator-bound}
		\|\mathcal K_\mu*f\|_{r'}\leq C\|f\|_r,
		\qquad
		|\mathfrak B_{\mathcal C}(f,g)|\leq C\|f\|_r\|g\|_r.
	\end{equation}
	Set $f=|w|^p$, $g_n=|z_n|^p$, $h_0=|w|^{p-2}w$, and $h_n=|z_n|^{p-2}z_n$.
	Applying the standard Br\'{e}zis--Lieb lemma \cite{BrezisLieb1983}, we obtain
	\begin{equation}\label{eq:BL-two-remainders}
		|w_n|^p-f-g_n\to0\text{ in }L^r(\mathcal C),
		\qquad
		|w_n|^{p-2}w_n-h_0-h_n\to0\text{ in }L^h(\mathcal C).
	\end{equation}
	
    Since $\{g_n\}$ is bounded in $L^r(\mathcal C)$ and $g_n\to0$ a.e., we have $g_n\rightharpoonup0$ in $L^r(\mathcal C)$. Combined with $\mathcal K_\mu*f\in L^{r'}(\mathcal C)$, this guarantees $\mathfrak B_{\mathcal C}(f,g_n)\to0$. The conclusion \eqref{eq:D-BL-value} then follows by expanding the bilinear form and applying the first statement in \eqref{eq:BL-two-remainders}.
	
	To prove the convergence of the derivatives, we estimate the two mixed terms. Let
	$K_L=[-L,L]\times\mathbb S^{N-1}$. Since $\ell=2^{\sharp}<2^*$, by using Rellich compactness on $K_L$, we obtain
	\begin{equation}\label{eq:BL-local-powers}
		|g_n|_{L^r(K_L)}\to0,\qquad
		|h_n|_{L^h(K_L)}\to0
	\end{equation}
	for every fixed $L>0$. As $\mathcal K_\mu*f$ is fixed, we apply H\"older's inequality to deduce
	\[
	\|(\mathcal K_\mu*f)h_n\|_{\ell'}
	\leq\|\mathcal K_\mu*f\|_{L^{r'}(K_L)}\|h_n\|_{L^h(K_L)}
	+C\|\mathcal K_\mu*f\|_{L^{r'}(\mathcal C\setminus K_L)}.
	\]
	Letting $n\to\infty$ and $L\to\infty$, we obtain
	\begin{equation}\label{eq:BL-first-mixed}
		\|(\mathcal K_\mu*f)h_n\|_{\ell'}\to0.
	\end{equation}
	
	We next show that $\mathcal K_\mu*g_n\to0$ in $L^{r'}(K_L)$. Fix $R>L+2$. The local part
	$\mathcal K_\mu*(g_n\mathbf1_{K_R})$ tends to zero in $L^{r'}(\mathcal C)$ by
	\eqref{eq:BL-operator-bound} and \eqref{eq:BL-local-powers}. For the tail part, taking $z=(t,\omega)\in K_L$, we apply H\"older's inequality to obtain
	\[
	\left|\bigl(\mathcal K_\mu*(g_n\mathbf1_{\mathcal C\setminus K_R})\bigr)(z)\right|
	\leq\|g_n\|_r
	\left(\int_{|s|>R}\int_{\mathbb S^{N-1}}
	\mathcal K_\mu(t-s,\omega,\omega')^{r'}\,d\omega'\,ds\right)^{1/r'}
	\leq C e^{-\mu(R-L)/2}.
	\]
	The last estimate uses $|t-s|>2$ and the exponential axial decay of $\mathcal K_\mu$, and it is uniform in $z\in K_L$ and $n$. Letting $n\to\infty$ and $R\to\infty$, we find
	\[
	\lim_{n\to\infty}\|\mathcal K_\mu*g_n\|_{L^{r'}(K_L)}=0.
	\]
	Since the sequence $\{\mathcal K_\mu*g_n\}$ is bounded in $L^{r'}(\mathcal C)$, we have
	\[
	\|(\mathcal K_\mu*g_n)h_0\|_{\ell'}
	\leq\|\mathcal K_\mu*g_n\|_{L^{r'}(K_L)}\|h_0\|_{L^h(K_L)}
	+C\|h_0\|_{L^h(\mathcal C\setminus K_L)}\longrightarrow0
	\]
	by letting $n\to\infty$ and $L\to\infty$. Thus
	\begin{equation}\label{eq:BL-second-mixed}
		\|(\mathcal K_\mu*g_n)h_0\|_{\ell'}\to0.
	\end{equation}
	Expanding $(\mathcal K_\mu*|w_n|^p)|w_n|^{p-2}w_n$ via
	\eqref{eq:BL-two-remainders}, every remainder tends to zero in $L^{\ell'}(\mathcal C)$ by H\"older's inequality and \eqref{eq:BL-operator-bound}. The remaining mixed terms vanish by
	\eqref{eq:BL-first-mixed} and \eqref{eq:BL-second-mixed}. The continuous embedding $H^1(\mathcal C)\hookrightarrow L^\ell(\mathcal C)$ yields the splitting property in $H^{-1}(\mathcal C)$. By duality, this property naturally carries over to $\bigl(D_a^{1,2}(\mathbb R^N)\bigr)^*$.
	
    Finally, for a fixed test function $\varphi\in C_c^\infty(\mathcal C)$, by combining the boundedness of $\mathcal K_\mu*g_n$ in $L^{r'}(\mathcal C)$ and the local strong convergence of $h_n$ in $L^h$ from \eqref{eq:BL-local-powers}, we deduce that
    $\int_{\mathcal C}(\mathcal K_\mu*g_n)h_n\varphi\,dt\,d\omega\to0$. For every $\varphi\in H^1(\mathcal C)$, we have
    \[
    \left|\mathcal D_{\mathcal C}'(z_n)[\varphi]\right|
    \leq C\|z_n\|_{H^1(\mathcal C)}^{2p-1}
    \|\varphi\|_{H^1(\mathcal C)}.
    \]
    Since $(z_n)$ is bounded in $H^1(\mathcal C)$,
    the functionals $\mathcal D_{\mathcal C}'(z_n)$ are
    uniformly bounded in $H^{-1}(\mathcal C)$.
    Density then gives
    $\mathcal D_{\mathcal C}'(z_n)[\varphi]\to0$
    for every fixed $\varphi\in H^1(\mathcal C)$.
\end{proof}

\subsection{Mass-preserving scaling and the Pohozaev identity}\label{sec:pohozaev}
Throughout this subsection, we assume \eqref{eq:decay-range},
which is equivalent to $p>2$. To derive the Pohozaev identity,
we use the cylindrical formulation. Let $u\in X$ be a weak
solution of the Euler--Lagrange equation with Lagrange
multiplier $\lambda\in\mathbb R$, and let $w$ be its
Emden--Fowler transform. The corresponding equation on
$\mathcal C$ satisfies
\begin{equation}\label{eq:normalized-cylinder-equation}
		-\partial_{tt}w-\Delta_{\omega}w+\frac{\kappa^2}{4}w
		=\lambda e^{-2t}w+\beta e^{-\gamma_qt}|w|^{q-2}w
		+(\mathcal K_\mu*|w|^p)|w|^{p-2}w
		\quad\text{in }\mathcal C.
\end{equation}

To establish the regularity and asymptotic decay needed to justify the Pohozaev identity, we first prove a uniform local $L^\infty$ estimate via Moser iteration.

\begin{lemma}[Local boundedness]\label{lem:local-moser}
	Let $K_R(t_0)=(t_0-R,t_0+R)\times\mathbb S^{N-1}$. Suppose that
	$z\in H^1(K_2(t_0))$ is nonnegative and satisfies $L_\kappa z\leq Qz$ weakly in $K_2(t_0)$,
	where $0\leq Q\in L^s(K_2(t_0))$ with $\|Q\|_{L^s(K_2(t_0))}\leq B$ for some $s>N/2$.
	Then there exists a constant $C=C(B,s,N)>0$, independent of $t_0$ and $\kappa$, such that
	\begin{equation}\label{eq:uniform-local-boundedness}
		\|z\|_{L^\infty(K_1(t_0))}
		\leq C(B,s,N)\|z\|_{L^2(K_2(t_0))}.
	\end{equation}
\end{lemma}
\begin{proof}
We write $\nabla_{\mathcal C}=(\partial_t,\nabla_\omega)$
for the cylindrical gradient. Fix $m\geq1$, $L>0$,
and $1\leq r<R\leq2$. Choose
$\xi\in C_c^\infty(K_R(t_0))$ such that
$0\leq\xi\leq1$, $\xi=1$ on $K_r(t_0)$, and
$
|\nabla_{\mathcal C}\xi|\leq\frac{C}{R-r}.
$
Set $z_L=\min\{z,L\}$ and $y=zz_L^{m-1}$.
Since $z\geq0$, the function
$\xi^2zz_L^{2m-2}$ is a nonnegative admissible
test function in $H_0^1(K_R(t_0))$. We test $L_\kappa z\leq Qz$ with
$\xi^2zz_L^{2m-2}$. The potential term
$\frac{\kappa^2}{4}z$ and the contribution
\[
(2m-2)\int_{K_R(t_0)\cap\{z<L\}}
\xi^2z^{2m-2}|\nabla_{\mathcal C}z|^2
\,dt\,d\omega
\]
are both nonnegative and can be discarded.
Applying Young's inequality to the cutoff term,
we obtain
	\[
	\int_{K_R(t_0)}\xi^2z_L^{2m-2}|\nabla_{\mathcal C}z|^2\,dt\,d\omega
	\leq 2\int_{K_R(t_0)} Q\xi^2y^2\,dt\,d\omega
	+4\int_{K_R(t_0)}|\nabla_{\mathcal C}\xi|^2y^2\,dt\,d\omega.
	\]
Since
$|\nabla_{\mathcal C}y|
\leq m z_L^{m-1}|\nabla_{\mathcal C}z|$
almost everywhere, the local Sobolev inequality gives
\[
\|\xi y\|_{2^*}^2
\leq C m^2
\left(
\int_{K_R(t_0)}Q(\xi y)^2\,dt\,d\omega
+
\int_{K_R(t_0)}
(|\nabla_{\mathcal C}\xi|^2+1)y^2\,dt\,d\omega
\right).
\]
To estimate the first integral, set
$\vartheta=N/(2s)\in(0,1)$ and $s'=s/(s-1)$.
By H\"older's inequality and interpolation, for any
$v\in H_0^1(K_R(t_0))$ we have
\[
\int_{K_R(t_0)}Qv^2\,dt\,d\omega
\leq B\|v\|_{2s'}^2
\leq B\|v\|_{2^*}^{2\vartheta}
\|v\|_2^{2(1-\vartheta)}.
\]
Young's inequality then yields, for every $\epsilon>0$,
\[
\int_{K_R(t_0)}Qv^2\,dt\,d\omega
\leq
\epsilon\|v\|_{2^*}^2
+
C_\vartheta B^{1/(1-\vartheta)}
\epsilon^{-\vartheta/(1-\vartheta)}
\|v\|_2^2.
\]
We apply this estimate with $v=\xi y$ and choose
$\epsilon=(2Cm^2)^{-1}$. Absorbing the corresponding
term into the left-hand side and using $R-r\leq1$,
we obtain
\[
\|\xi y\|_{2^*}^2
\leq C_B m^{b_0}(R-r)^{-2}
\|y\|_{L^2(K_R(t_0))}^2,
\qquad
b_0=\frac{2}{1-\vartheta}.
\]

Suppose that $z\in L^{2m}(K_R(t_0))$.
We may let $L\to\infty$ in the preceding estimate.
Since $\xi=1$ on $K_r(t_0)$, monotone convergence gives
\[
\|z\|_{L^{2m\chi}(K_r(t_0))}
\leq
\left[C_Bm^{b_0}(R-r)^{-2}\right]^{1/(2m)}
\|z\|_{L^{2m}(K_R(t_0))},
\qquad
\chi=\frac{N}{N-2}.
\]

We now iterate this inequality with
$m_j=\chi^j$ and $R_j=1+2^{-j}$, starting from
$m_0=1$. After $j+1$ steps, we obtain
\[
\|z\|_{L^{2\chi^{j+1}}(K_{R_{j+1}}(t_0))}
\leq
\prod_{i=0}^{j}
\left[
C_B\chi^{ib_0}2^{2i+2}
\right]^{1/(2\chi^i)}
\|z\|_{L^2(K_2(t_0))}.
\]
The product is bounded independently of $j$, since
\[
\sum_{i=0}^{\infty}(1+i)\chi^{-i}<\infty.
\]
Letting $j\to\infty$, we obtain
\eqref{eq:uniform-local-boundedness}.
The local Sobolev constant is independent of $t_0$
by the translation invariance of the cylinder
$\mathcal C$.
\end{proof}

\begin{lemma}[Local regularity and right tail]\label{lem:full-equation-regularity}
	Every real weak solution $u\in X$ of the Euler--Lagrange equation in \eqref{eq:problem} has Emden--Fowler transform
	$w\in W^{2,r}_{\mathrm{loc}}(\mathcal C)$ for all $r\in[1,\infty)$, and hence $w\in C^{1,\sigma}_{\mathrm{loc}}(\mathcal C)$ for some $\sigma\in(0,1)$. Moreover, both $w$ and $(\mathcal K_\mu*|w|^p)$ are uniformly bounded on $[0,\infty)\times\mathbb S^{N-1}$, and
	\begin{equation}\label{eq:normalized-right-vanishing}
		\sup_{\omega\in\mathbb S^{N-1}}|w(t,\omega)|\longrightarrow0,\quad\text{ }t\to+\infty.
	\end{equation}
\end{lemma}
\begin{proof}
	Set $V=\mathcal K_\mu*|w|^p$. For any nonnegative
	$\varphi\in C_c^\infty(\mathcal C)$, we test
	\eqref{eq:normalized-cylinder-equation} with
	$\varphi w/(w^2+\varepsilon^2)^{1/2}$ and let
	$\varepsilon\to0^+$. Kato's inequality gives
	$L_\kappa|w|\leq Q|w|$ in $\mathcal D'(\mathcal C)$, where
	\[
	Q=|\lambda|e^{-2t}
	+\beta e^{-\gamma_qt}|w|^{q-2}
	+V|w|^{p-2}.
	\]
	
	Set $s_0=2^\sharp/(2^\sharp-2)=N/(2d)>N/2$.
	Since $1/s_0=\mu/(2N)+(p-2)/2^\sharp$ and
	$\|V\|_{2N/\mu}\leq C\|w\|_{2^\sharp}^p$,
	H\"older's inequality gives
	$\|V|w|^{p-2}\|_{s_0}\leq C\|w\|_{2^\sharp}^{2p-2}$.
	Moreover, $2<q<2^\sharp$ implies
	$2^\sharp/(q-2)>s_0$, and hence
	\[
	\||w|^{q-2}\|_{L^{s_0}(K_2(t_0))}
	\leq C\|w\|_{2^\sharp}^{q-2}
	\]
	uniformly in $t_0\in\mathbb R$.
	Thus $Q\in L^{s_0}_{\mathrm{loc}}(\mathcal C)$.
	Since the exponential weights are uniformly bounded
	on $K_2(t_0)$ for $t_0\geq0$, we also have
	$\sup_{t_0\geq0}\|Q\|_{L^{s_0}(K_2(t_0))}<\infty$.
	Applying \eqref{eq:uniform-local-boundedness}
	to $|w|$, we obtain $w\in L^\infty_{\mathrm{loc}}(\mathcal C)$
	and
	\[
	\|w\|_{L^\infty(K_1(t_0))}
	\leq C\|w\|_{L^2(K_2(t_0))},
	\qquad t_0\geq0,
	\]
	with $C$ independent of $t_0$.
	Since $w\in L^2(\mathcal C)$, the right-hand side
	tends to zero as $t_0\to+\infty$.
	This proves \eqref{eq:normalized-right-vanishing}
	and the uniform boundedness of $w$ on
	$[0,\infty)\times\mathbb S^{N-1}$.
	
	We next estimate $V$. For $t\geq0$, we split
	its defining integral into $|t-s|\leq1$ and
	$|t-s|>1$. In the first region, $s\geq-1$,
	where $w$ is uniformly bounded. The local
	integrability of $\mathcal K_\mu$ gives a
	uniform bound for this contribution.
	In the second region, H\"older's inequality yields
	\[
	\int_{|t-s|>1}\int_{\mathbb S^{N-1}}
	\mathcal K_\mu(t-s,\omega,\omega')
	|w(s,\omega')|^p\,d\omega'\,ds
	\leq C\||w|^p\|_{L^{r_\mu}(\mathcal C)},
	\]
	since the kernel has a uniformly bounded
	$L^{r_\mu'}$ norm away from $t=s$.
	Hence $V$ is uniformly bounded on
	$[0,\infty)\times\mathbb S^{N-1}$.
	The same argument on bounded slabs gives
	$V\in L^\infty_{\mathrm{loc}}(\mathcal C)$.
	
	The right-hand side of
	\eqref{eq:normalized-cylinder-equation}
	therefore belongs to $L^\infty_{\mathrm{loc}}(\mathcal C)$.
	Standard interior elliptic estimates give
	$w\in W^{2,r}_{\mathrm{loc}}(\mathcal C)$
	for every $1\leq r<\infty$.
	Taking $r>N$, the Sobolev embedding yields
	$w\in C^{1,\sigma}_{\mathrm{loc}}(\mathcal C)$,
	where $\sigma=1-N/r$.
\end{proof}
\begin{proposition}[Pohozaev identity]\label{prop:pohozaev}
	Every weak solution $u\in X$ satisfies
	\begin{equation}\label{eq:lambda-pohozaev-cylinder}
		\lambda M(u)=\beta(\delta_q-1)B_q(u),
	\end{equation}
	and
	\[
	P_{\beta}(u)=0.
	\]
	Consequently, every nonzero normalized solution has $\lambda<0$.
\end{proposition}

\begin{proof}
	Let $\eta\in C_c^\infty(\mathbb R)$ satisfy
	$0\leq\eta\leq1$, $\eta=1$ on $[-1,1]$, and
	$\eta=0$ outside $[-2,2]$. For $R\geq1$, set
	$\eta_R(t)=\eta(t/R)$. By Lemma~\ref{lem:full-equation-regularity},
	$\eta_R\partial_tw\in H^1(\mathcal C)$ has compact
	axial support, and its inverse Emden--Fowler transform
	belongs to $X$. We may therefore use
	$\eta_R\partial_tw$ as a test function in the
	cylindrical weak formulation.
	
	We integrate by parts to obtain
	\[
	\int_{\mathcal C}L_\kappa w\,\eta_R\partial_tw
	\,dt\,d\omega
	=\frac12\int_{\mathcal C}\eta_R'
	\left(|\partial_tw|^2-|\nabla_\omega w|^2
	-\frac{\kappa^2}{4}w^2\right)\,dt\,d\omega.
	\]
	This tends to zero as $R\to\infty$, since
	$\|\eta_R'\|_\infty\leq C/R$ and $w\in H^1(\mathcal C)$.
	For the weighted linear and local terms, we obtain
	\[
	\begin{aligned}
		\lambda\int_{\mathcal C}e^{-2t}w\eta_R\partial_tw\,dt\,d\omega
		&=\lambda\int_{\mathcal C}e^{-2t}\eta_Rw^2\,dt\,d\omega
		-\frac{\lambda}{2}\int_{\mathcal C}e^{-2t}\eta_R'w^2\,dt\,d\omega,\\
		\beta\int_{\mathcal C}e^{-\gamma_qt}|w|^{q-2}w
		\eta_R\partial_tw\,dt\,d\omega
		&=\frac{\beta\gamma_q}{q}
		\int_{\mathcal C}e^{-\gamma_qt}\eta_R|w|^q\,dt\,d\omega
		-\frac{\beta}{q}
		\int_{\mathcal C}e^{-\gamma_qt}\eta_R'|w|^q\,dt\,d\omega,
	\end{aligned}
	\]
	By \eqref{eq:cylinder-mass-local}, the terms
	involving $\eta_R'$ vanish as $R\to\infty$,
	while the remaining terms converge to
	$\lambda M(u)$ and $(\beta\gamma_q/q)B_q(u)$.
	
	It remains to treat the nonlocal term. Set
	$f=|w|^p$ and
	\[
	I_R=\frac1p\int_{\mathcal C}
	(\mathcal K_\mu*f)\eta_R\partial_tf\,dt\,d\omega.
	\]
	For $0<\varepsilon\leq1$, introduce
	\[
	\mathcal K_{\mu,\varepsilon}(t-s,\omega,\omega')
	=\left[2(\cosh(t-s)-\omega\cdot\omega')
	+\varepsilon^2\right]^{-\mu/2},
	\]
	and let $I_{R,\varepsilon}$ be defined as $I_R$
	with $\mathcal K_\mu$ replaced by
	$\mathcal K_{\mu,\varepsilon}$.
	In the following double integrals, kernels without
	displayed arguments are evaluated at
	$(t-s,\omega,\omega')$.
	
	For fixed $R$ and $\varepsilon$, set	
	\[
	I_{R,\varepsilon}
	=\frac1p\int_{\mathcal C}(\mathcal K_{\mu,\varepsilon}*f)
	\eta_R\partial_tf\,dt\,d\omega.
	\]
	We integrate by parts and then symmetrize in $(t,\omega)$ and $(s,\omega')$ to obtain
	\[
		I_{R,\varepsilon}
		=-\frac1{2p}\iint_{\mathcal C\times\mathcal C}
		\Bigl[
		\mathcal K_{\mu,\varepsilon}
		\bigl(\eta_R'(t)+\eta_R'(s)\bigr)
		+\partial_t\mathcal K_{\mu,\varepsilon}
		\bigl(\eta_R(t)-\eta_R(s)\bigr)
		\Bigr]f(t,\omega)f(s,\omega')
		\,dt\,d\omega\,ds\,d\omega'.
	\]
	Here we have used
	$\partial_s\mathcal K_{\mu,\varepsilon}
	=-\partial_t\mathcal K_{\mu,\varepsilon}$.
	The auxiliary cutoff can be removed by the
	cylindrical HLS inequality, since
	$\mathcal K_{\mu,\varepsilon}\leq\mathcal K_\mu$
	and $|\partial_t\mathcal K_{\mu,\varepsilon}|
	\leq C_\varepsilon\mathcal K_\mu$.
	
	To pass to the limit as $\varepsilon\to0^+$, we denote by
	$d_{\mathcal C}$ the distance between $(t,\omega)$ and
	$(s,\omega')$ on $\mathcal C$. Using
	\[
	\partial_t\mathcal K_{\mu,\varepsilon}
	=-\mu\frac{\sinh(t-s)}
	{2(\cosh(t-s)-\omega\cdot\omega')+\varepsilon^2}
	\mathcal K_{\mu,\varepsilon}
	\]
	and the estimate
	$|\eta_R(t)-\eta_R(s)|\leq C|t-s|/R$, we obtain
	\[
	|\eta_R(t)-\eta_R(s)|
	|\partial_t\mathcal K_{\mu,\varepsilon}|
	\leq\frac{C}{R}
	\frac{|t-s|^2}{d_{\mathcal C}^2+\varepsilon^2}
	\mathcal K_{\mu,\varepsilon}
	\leq\frac{C}{R}\mathcal K_\mu
	\quad\text{if }d_{\mathcal C}\leq1.
	\]
	On $\{d_{\mathcal C}>1\}$, we have
	$|\partial_t\mathcal K_{\mu,\varepsilon}|
	\leq C\mathcal K_\mu$, uniformly for
	$0<\varepsilon\leq1$.
	These estimates and the cylindrical HLS inequality
	provide an integrable majorant for the symmetrized
	integral. Since $\eta_R\partial_tf$ is bounded
	and compactly supported, we also have
	$I_{R,\varepsilon}\to I_R$.
	We apply the dominated convergence theorem to obtain
	\begin{equation}\label{eq:nonlocal-translation-cutoff}
			I_R=-\frac1{2p}\iint_{\mathcal C\times\mathcal C}
			\Bigl[
			\mathcal K_\mu
			\bigl(\eta_R'(t)+\eta_R'(s)\bigr)
			+\partial_t\mathcal K_\mu
			\bigl(\eta_R(t)-\eta_R(s)\bigr)
			\Bigr]f(t,\omega)f(s,\omega')
			\,dt\,d\omega\,ds\,d\omega'.
	\end{equation}
	
	We now let $R\to\infty$. By the cylindrical
	HLS inequality, the terms involving $\eta_R'$
	are bounded by $CR^{-1}\|f\|_{r_\mu}^2$.
	The same bound holds for the kernel derivative
	term on $\{d_{\mathcal C}\leq1\}$.
	On $\{d_{\mathcal C}>1\}$, we use
	$|\partial_t\mathcal K_\mu|\leq C\mathcal K_\mu$
	and
	\[
	|\eta_R(t)-\eta_R(s)|
	\leq\mathbf1_{\{|t|>R\}}
	+\mathbf1_{\{|s|>R\}}.
	\]
	Writing $K_R=[-R,R]\times\mathbb S^{N-1}$,
	another application of the cylindrical HLS
	inequality gives
	\[
	|I_R|\leq\frac{C}{R}\|f\|_{r_\mu}^2
	+C\|f\mathbf1_{\mathcal C\setminus K_R}\|_{r_\mu}
	\|f\|_{r_\mu}\longrightarrow0.
	\]
	
	Passing to the limit in the weak formulation,
	we obtain
	\[
	\lambda M(u)+\frac{\beta\gamma_q}{q}B_q(u)=0.
	\]
	Since $\gamma_q/q=1-\delta_q$, this proves
	\eqref{eq:lambda-pohozaev-cylinder}.
	Testing the original equation with $u$, we also have
	\begin{equation}\label{eq:nehari-identity-solution}
		A(u)=\lambda M(u)+\beta B_q(u)+\mathcal D(u).
	\end{equation}
	Combining the two identities yields $P_\beta(u)=0$.
	Finally, $\beta>0$, $\delta_q\in(0,1)$, and
	$M(u),B_q(u)>0$ for every nontrivial solution,
	so \eqref{eq:lambda-pohozaev-cylinder} implies
	$\lambda<0$.
\end{proof}

\begin{remark}\label{rem:classical-pohozaev-equivalence}
	Combining the cylindrical Pohozaev identity \eqref{eq:lambda-pohozaev-cylinder} with the Nehari identity \eqref{eq:nehari-identity-solution} recovers the Euclidean Pohozaev identity
	\[
	\frac{N-2-2a}{2}A(u)
	=\frac{N-2a}{2}\lambda M(u)
	+\frac{N-bq}{q}\beta B_q(u)
	+\frac{2N-2bp-\mu}{2p}\mathcal D(u).
	\]
	In cylindrical coordinates, this identity is obtained by using compactly supported axial derivatives as test functions, where the nonlocal commutator is controlled via \eqref{eq:nonlocal-translation-cutoff}.
\end{remark}

\subsection{Fiber geometry and the Pohozaev manifold}\label{sec:fibers}

Throughout this subsection assume \eqref{eq:decay-range}, $\rho,\beta>0$, and $2<q<2^\sharp$. On $S_{\rho}$, we have
\begin{equation}\label{eq:Bq-shell}
	B_q(u)\leq K_qA(u)^{r_q/2},
	\qquad K_q=C_{a,b,q}\rho^{\gamma_q}.
\end{equation}

\begin{lemma}[Fiber maps in the mass-subcritical case]\label{lem:subcritical-fibers}
	Assume $2<q<q_c$.  There exists $\beta_{\mathrm{fib}}>0$ such that, for $0<\beta<\beta_{\mathrm{fib}}$ and every $u\in S_{\rho}$, $\Phi_{\beta,u}$ has exactly two critical points
	\[
	t_1(u)<t_2(u).
	\]
	The maps $u\mapsto t_i(u)$ are of class $C^1$, and
	\begin{equation}\label{eq:branch-projections}
		t_1(u)\star u\in\mathcal M_{\rho,\beta}^{+},
		\qquad
		t_2(u)\star u\in\mathcal M_{\rho,\beta}^{-}.
	\end{equation}
	The first critical point is a strict local minimum of the fiber map at a negative level, whereas the second is its strict global maximum at a positive level. Moreover, $\mathcal M_{\rho,\beta}^{0}=\varnothing$.
\end{lemma}

\begin{proof}
	Set $s=e^t$ and define
	$f_u(s)=A(u)s^{2-r_q}-\mathcal D(u)s^{2p-r_q}$ for $s>0$.
	The equation $\Phi_{\beta,u}'(t)=0$ is equivalent to
	\begin{equation}\label{eq:fiber-equation-s}
		f_u(s)=\beta\delta_qB_q(u).
	\end{equation}
	Since $0<r_q<2<2p$, we have $f_u(0)=0$ and
	$f_u(s)\to-\infty$ as $s\to\infty$. Moreover,
	\begin{equation}\label{eq:fiber-scalar-derivative-proof}
		f_u'(s)=s^{1-r_q}
		\left((2-r_q)A(u)
		-(2p-r_q)\mathcal D(u)s^{2p-2}\right).
	\end{equation}
	Thus $f_u$ has a unique critical point, which is
	its strict maximum.
	
	Set $z=sA(u)^{1/2}$. By the critical Stein--Weiss
	inequality, $\mathcal D(u)\leq S_\mu(a,b)^{-p}A(u)^p$,
	and hence
	$f_u(s)\geq A(u)^{r_q/2}g(z)$, where
	$g(z)=z^{2-r_q}-S_\mu(a,b)^{-p}z^{2p-r_q}$.
	Let $g_*=\max_{z>0}g(z)>0$. On the other hand,
	\eqref{eq:Bq-shell} gives
	$\beta\delta_qB_q(u)\leq
	\beta\delta_qK_qA(u)^{r_q/2}$.
	Choose $\beta_{\mathrm{fib}}>0$ such that
	$\beta_{\mathrm{fib}}\delta_qK_q<g_*$.
	Then
	$\max_{s>0}f_u(s)>\beta\delta_qB_q(u)$
	for every $u\in S_\rho$ and
	$0<\beta\leq\beta_{\mathrm{fib}}$.
	Since $f_u$ increases before its maximum and
	decreases afterwards, \eqref{eq:fiber-equation-s}
	has exactly two roots $s_1(u)<s_2(u)$.
	These give two distinct critical points
	$t_i(u)=\log s_i(u)$.
	
	As $t\to-\infty$, the leading term of
	$\Phi_{\beta,u}(t)$ is
	$-(\beta/q)e^{r_qt}B_q(u)$.
	Hence $t_1(u)$ is a strict local minimum with
	negative energy. To obtain a positive lower bound
	for the second critical value, the same estimates give
	\begin{equation}\label{eq:fiber-positive-test-proof}
		\Phi_{\beta,u}(t)
		\geq \frac{z^2}{2}
		-\frac{\beta K_q}{q}z^{r_q}
		-\frac{S_\mu(a,b)^{-p}}{2p}z^{2p},
		\qquad z=e^tA(u)^{1/2}.
	\end{equation}
	Choose $z_0>0$ sufficiently small that
	$z_0^2/2-S_\mu(a,b)^{-p}z_0^{2p}/(2p)>z_0^2/4$.
	After decreasing $\beta_{\mathrm{fib}}$, we may also
	assume $\beta K_qz_0^{r_q}/q<z_0^2/8$.
	It follows from \eqref{eq:fiber-positive-test-proof}
	that $\Phi_{\beta,u}(t)>z_0^2/8$ whenever
	$e^tA(u)^{1/2}=z_0$.
	This lower bound is independent of $u$.
	Since $\Phi_{\beta,u}(t)\to0^-$ as $t\to-\infty$
	and $\Phi_{\beta,u}(t)\to-\infty$ as $t\to+\infty$,
	the second critical point is a strict global maximum
	with uniformly positive energy.
	
	At either root, the sign of $f_u'(s_i(u))$
	agrees with that of $Q_\beta(t_i(u)\star u)$;
	indeed,
	$Q_\beta(t_i(u)\star u)
	=s_i(u)^{r_q+1}f_u'(s_i(u))$.
	This proves \eqref{eq:branch-projections} and shows
	that $\mathcal M_{\rho,\beta}^{0}=\varnothing$.
	Finally, the derivative of
	$\Phi_{\beta,u}'(t)$ with respect to $t$ is nonzero
	at either critical point. The implicit function
	theorem on $S_\rho$ therefore gives the $C^1$
	dependence of $t_1$ and $t_2$ on $u$.
\end{proof}

\begin{lemma}[Fiber maps in the mass-critical and mass-supercritical cases]\label{lem:critical-super-fibers}
	Assume $q_c\leq q<2^{\sharp}$.  If $q=q_c$, assume $0<\beta<\beta_*$.  Then every fiber has exactly one critical point $t(u)$, which is its strict global maximum.  The map $u\mapsto t(u)$ is of class $C^1$, and
	\[
	Q_{\beta}(v)<0\quad\text{for every }v\in\mathcal M_{\rho,\beta}.
	\]
\end{lemma}

\begin{proof}
	If $r_q>2$, the critical-point equation becomes
	\[
	A(u)=\beta\delta_q e^{(r_q-2)t}B_q(u)
	+e^{(2p-2)t}\mathcal D(u).
	\]
	The right-hand side is strictly increasing from
	zero to infinity, so the equation has a unique solution
	$t(u)$. The sign of $\Phi_{\beta,u}'$ then shows that
	$t(u)$ is the unique strict global maximizer of
	$\Phi_{\beta,u}$.
	
	If $r_q=2$, then $q=q_c$ and $\delta_{q_c}=2/q_c$.
	By the definition of $\mathcal C_{q_c}$ in
	\eqref{C_{q_c}}, we have
	\begin{equation}\label{eq:critical-fiber-positive-coefficient}
	\frac12A(u)-\frac{\beta}{q_c}B_{q_c}(u)
	\geq
	\left(\frac12-
	\frac{\beta\mathcal C_{q_c}\rho^{q_c-2}}{q_c}\right)
	A(u)>0.
	\end{equation}
	Hence $\Phi_{\beta,u}(t)$ has the form
	$ce^{2t}-\mathcal D(u)e^{2pt}/(2p)$ with $c>0$,
	and therefore has a unique strict global maximum.
	
	For any $v\in\mathcal M_{\rho,\beta}$, we use \eqref{eq:Q-definition}
	to obtain
	\[
	Q_\beta(v)=
	\begin{cases}
		-(r_q-2)\beta\delta_qB_q(v)
		-(2p-2)\mathcal D(v),&r_q>2,\\
		-(2p-2)\mathcal D(v),&r_q=2.
	\end{cases}
	\]
	Thus $Q_\beta(v)<0$ in both cases. The implicit
	function theorem consequently yields
	$u\mapsto t(u)\in C^1(S_\rho)$.
\end{proof}
\begin{lemma}[Regularity and reduced functionals]\label{lem:natural-constraint}
	Under the assumptions of Lemma \ref{lem:subcritical-fibers} or Lemma \ref{lem:critical-super-fibers}, the set $\mathcal M_{\rho,\beta}$ is a $C^1$ submanifold of codimension one in $S_\rho$. In the mass-subcritical case, its two branches are open submanifolds. Let $t_i$ be one of the branch projections in Lemma \ref{lem:subcritical-fibers}, or let $t$ be the projection in Lemma \ref{lem:critical-super-fibers}.  Define
	\[
	\Psi_i(v)=E_{\beta}(t_i(v)\star v),
	\qquad
	\Psi(v)=E_{\beta}(t(v)\star v).
	\]
	If $v$ is a critical point of the corresponding reduced functional on $S_{\rho}$, then the projected point is a critical point of $E_{\beta}|_{S_{\rho}}$.
\end{lemma}

\begin{proof}
	The functional $P_\beta$ is $C^1$ on $X$. Suppose that
	$u\in\mathcal M_{\rho,\beta}$ and
	$d(P_\beta|_{S_\rho})(u)=0$. Since $M'(u)\neq0$,
	there exists $\tau\in\mathbb R$ such that
	$P_\beta'(u)=\tau M'(u)$. Its Emden--Fowler transform
	$w$ then satisfies
	\[
	L_\kappa w=\tau e^{-2t}w
	+\frac{\beta r_q}{2}e^{-\gamma_qt}|w|^{q-2}w
	+p(\mathcal K_\mu*|w|^p)|w|^{p-2}w.
	\]
	The regularity argument in Lemma
	\ref{lem:full-equation-regularity} applies with
	coefficients $\beta r_q/2$ and $p$. In particular,
	$w\in H^2_{\mathrm{loc}}(\mathcal C)$, so the cutoff
	test $\eta_R\partial_tw$ used in Proposition
	\ref{prop:pohozaev} is admissible. Repeating that
	argument, with the nonlocal term multiplied by $p$,
	yields the same vanishing commutator estimate
	\eqref{eq:nonlocal-translation-cutoff}. Combining the Pohozaev identity with the original equation tested against $u$, we obtain
	$\tau M(u)=-\beta r_q\gamma_q B_q(u)/(2q)$ and
	$A(u)=\tau M(u)+\beta r_q B_q(u)/2+p\mathcal D(u)$.
	Using $1-\gamma_q/q=\delta_q$ and $r_q=q\delta_q$,
	we eliminate $\tau M(u)$ to obtain
	\[
	Q_\beta(u)
	=2A(u)-\beta q\delta_q^2B_q(u)-2p\mathcal D(u)=0.
	\]
	This contradicts $Q_\beta(u)\neq0$, established in
	Lemma \ref{lem:subcritical-fibers} or
	Lemma \ref{lem:critical-super-fibers}. Hence zero is
	a regular value of $P_\beta|_{S_\rho}$ and
	$\mathcal M_{\rho,\beta}$ is a $C^1$ submanifold
	of $S_\rho$. Its two branches are open because
	$Q_\beta$ is continuous.
	
	To prove the second assertion, define
	\[
	F(t,v)=E_\beta(t\star v)
	=\frac{e^{2t}}2A(v)
	-\frac{\beta e^{r_qt}}qB_q(v)
	-\frac{e^{2pt}}{2p}\mathcal D(v).
	\]
	The explicit formula shows that
	$F\in C^1(\mathbb R\times S_\rho)$, without requiring
	differentiability of $t\mapsto t\star v$ in $X$.
	Moreover, $\partial_tF(t,v)=P_\beta(t\star v)$.
	Since $\Psi_i(v)=E_\beta(t_i(v)\star v)$, differentiating
	in the direction $\xi\in T_vS_\rho$ gives
	\begin{equation}\label{eq:reduced-derivative-formula}
		d\Psi_i(v)[\xi]
		=d(E_\beta|_{S_\rho})(u)[t_i(v)\star\xi]
		+t_i'(v)[\xi]P_\beta(u),
	\end{equation}
	where $u=t_i(v)\star v$.
	Since $u\in\mathcal M_{\rho,\beta}$, the last term
	vanishes. Furthermore, $\xi\mapsto t_i(v)\star\xi$
	is a bounded linear isomorphism from $T_vS_\rho$
	onto $T_uS_\rho$. Thus $d\Psi_i(v)=0$ implies
	$d(E_\beta|_{S_\rho})(u)=0$. The same argument
	applies to $\Psi$.
\end{proof}

\begin{lemma}[Uniform projection estimates]\label{lem:uniform-projection}
	Let $\widehat u_n$ lie on the upper branch, or on the unique branch when $q\geq q_c$. Assume
	\begin{equation}\label{eq:projection-input-bounds}
		0<c\leq A(\widehat u_n)\leq C,\qquad Q_\beta(\widehat u_n)\leq-c.
	\end{equation}
	Set $s_n=-\tfrac12\log A(\widehat u_n)$ and $\widehat v_n=s_n\star\widehat u_n$.
	If $v_n\in S_\rho$ and $\|v_n-\widehat v_n\|_X\to0$, then the corresponding projection parameters satisfy
	\begin{equation}\label{eq:uniform-projection-conclusion}
		|t_i(v_n)+s_n|\leq C\|v_n-\widehat v_n\|_X,\qquad
		|t_i(v_n)|\leq C,\qquad
		Q_\beta(t_i(v_n)\star v_n)\leq-c/2
	\end{equation}
	for large $n$, after decreasing $c$ if necessary. The same statements hold for the unique projection $t$.
\end{lemma}
\begin{proof}
	By the scaling identities, $A(\widehat v_n)=1$,
	$|s_n|\leq C$, and $t_i(\widehat v_n)=-s_n$.
	Moreover, $v_n$ and $\widehat v_n$ remain in a bounded
	subset of $X$. Since $A$, $B_q$, and $\mathcal D$ are
	continuously differentiable with bounded derivatives
	on bounded subsets of $X$, we have
	\[
	|A(v_n)-A(\widehat v_n)|
	+|B_q(v_n)-B_q(\widehat v_n)|
	+|\mathcal D(v_n)-\mathcal D(\widehat v_n)|
	\leq C\|v_n-\widehat v_n\|_X.
	\]
	
	For $v\in S_\rho$, define
	\[
	F_v(t)=e^{2t}A(v)-\beta\delta_qe^{r_qt}B_q(v)
	-e^{2pt}\mathcal D(v).
	\]
	Since $\{s_n\}$ is bounded, we may choose a fixed
	compact interval $I$ containing a neighborhood of
	every $-s_n$. The preceding estimates then imply
	\[
	\|F_{v_n}-F_{\widehat v_n}\|_{C^1(I)}
	\leq C\|v_n-\widehat v_n\|_X.
	\]
	By construction,
	$F_{\widehat v_n}(-s_n)=0$ and
	$F_{\widehat v_n}'(-s_n)=Q_\beta(\widehat u_n)\leq-c$.
The uniform bound on $F_{\widehat v_n}''$ on $I$ allows us to choose $\eta>0$ independently of $n$ such that
	\[
	F_{\widehat v_n}'(t)\leq-\frac{3c}{4}
	\qquad\text{whenever }|t+s_n|\leq\eta.
	\]
	For large $n$, the $C^1$ estimate gives
	$F_{v_n}'(t)\leq-c/2$ throughout the same interval.
	Moreover, $F_{v_n}(-s_n-\eta)>0$ and
	$F_{v_n}(-s_n+\eta)<0$. Thus $F_{v_n}$ has a unique
	zero in this interval, corresponding to the upper
	branch. We apply the mean value theorem to obtain
	\[
		|t_i(v_n)+s_n|
		\leq \frac{2}{c}|F_{v_n}(-s_n)|
		\leq C\|v_n-\widehat v_n\|_X.
	\]
	Since $\{s_n\}$ is bounded, $|t_i(v_n)|\leq C$.
	Furthermore,
	\[
	Q_\beta(t_i(v_n)\star v_n)
	=F_{v_n}'(t_i(v_n))\leq-\frac c2.
	\]
	The same argument applies to the unique projection
	when $q\geq q_c$.
\end{proof}
For later use, recall that the scaling identities give
\begin{equation}\label{eq:dilation-operator-norms}
	\begin{aligned}
		\|t\star\xi\|_X^2
		&=e^{2t}A(\xi)+M(\xi),\\
		\|t\star\xi\|_X
		&\leq\max\{1,e^t\}\|\xi\|_X,\\
		\|(-t)\star\xi\|_X
		&\leq\max\{1,e^{-t}\}\|\xi\|_X,
		\qquad \xi\in X.
	\end{aligned}
\end{equation}
These estimates also apply to the induced maps
between the tangent spaces $T_vS_\rho$ and
$T_{t\star v}S_\rho$.

\section{The critical problem and global compactness}\label{sec:critical-compactness}

In this section, we assume \eqref{eq:assumption-ab} and $0<\mu<N$. The additional restriction $p>2$ is imposed only for the decay estimates and the constrained problem.

\subsection{Existence and decay of critical extremals}\label{sec:extremals}

We begin by recalling a Lions-type vanishing lemma on the cylinder.

\begin{lemma}[Lions lemma on the cylinder]\label{lem:lions-cylinder}
	Let $\{w_n\}$ be bounded in $H^1(\mathcal C)$.  If
	\begin{equation}\label{eq:cylinder-vanishing}
		\sup_{s\in\mathbb R}
		\int_{(s-1,s+1)\times\mathbb S^{N-1}}|w_n|^2\,dt\,d\omega\to0,
	\end{equation}
	then
	\[
	w_n\to0\quad\text{in }L^r(\mathcal C)
	\]
	for every $2<r<2^*$.
\end{lemma}

\begin{proof}
	Set $r_0=2+4/N$. By interpolation and the uniform
	local Sobolev inequality, we have
	\[
	\int_{I_j}|w_n|^{r_0}\,dt\,d\omega
	\leq C
	\left(\int_{(j-1,j+2)\times\mathbb S^{N-1}}
	|w_n|^2\,dt\,d\omega\right)^{2/N}
	\|w_n\|_{H^1((j-1,j+2)\times\mathbb S^{N-1})}^2,
	\]
	where $I_j=(j,j+1)\times\mathbb S^{N-1}$.
	Each enlarged slab is covered by two slabs of the
	form appearing in \eqref{eq:cylinder-vanishing}.
	Hence the $L^2$ integral on the right tends to zero
	uniformly in $j$. Summing over $j\in\mathbb Z$ and
	using the boundedness of $\{w_n\}$ in $H^1(\mathcal C)$
	and the finite overlap of the enlarged slabs, we obtain
	$\|w_n\|_{L^{r_0}(\mathcal C)}\to0$. Interpolation with the bounded $L^2(\mathcal C)$ and $L^{2^*}(\mathcal C)$ norms then gives
	$w_n\to0$ in $L^r(\mathcal C)$ for every $2<r<2^*$.
\end{proof}

\begin{proposition}[Attainment of the critical quotient]\label{prop:attainment}
	The infimum $S_{\mu}(a,b)$ is attained.
\end{proposition}

\begin{proof}
	In cylindrical coordinates, write
	\[
	\|w\|_\kappa^2
	=\int_{\mathcal C}\left(
	|\partial_tw|^2+|\nabla_\omega w|^2
	+\frac{\kappa^2}{4}w^2\right)\,dt\,d\omega.
	\]
	Choose a minimizing sequence $(w_n)$ normalized by
	\begin{equation}\label{eq:minimizing-normalization}
		\mathcal D_{\mathcal C}(w_n)=1,\qquad
		\|w_n\|_\kappa^2\to S_\mu(a,b).
	\end{equation}
	By \eqref{eq:D-cylinder-bound}, the sequence is bounded
	away from zero in $L^{2^\sharp}(\mathcal C)$.
	The vanishing alternative is impossible. Lemma
	\ref{lem:lions-cylinder} gives $w_n\to0$ in
	$L^{2^\sharp}(\mathcal C)$ under this assumption.
	By \eqref{eq:D-cylinder-bound}, we then have
	$\mathcal D_{\mathcal C}(w_n)\to0$, contradicting
	$\mathcal D_{\mathcal C}(w_n)=1$.
	Hence there exist translations $s_n\in\mathbb R$ such that
	\begin{equation}\label{eq:nonvanishing-slab}
		\int_{(-1,1)\times\mathbb S^{N-1}}
		|w_n(t+s_n,\omega)|^2\,dt\,d\omega
		\geq c_0>0.
	\end{equation}
	Both $\|\cdot\|_\kappa$ and $\mathcal D_{\mathcal C}$
	are invariant under axial translations. Replacing $w_n$
	by $w_n(\cdot+s_n,\cdot)$ and passing to a subsequence,
	we may assume that $w_n\rightharpoonup w_*$ in
	$H^1(\mathcal C)$. Local compactness and
	\eqref{eq:nonvanishing-slab} ensure that
	$w_*\not\equiv0$.
	
	Set $z_n=w_n-w_*$. By weak convergence and the
	nonlocal Brezis--Lieb lemma, we have
	\begin{equation}\label{eq:minimizer-splitting}
		\begin{aligned}
			\|w_n\|_\kappa^2
			&=\|w_*\|_\kappa^2+\|z_n\|_\kappa^2+o(1),\\
			1&=\mathcal D_{\mathcal C}(w_*)
			+\mathcal D_{\mathcal C}(z_n)+o(1).
		\end{aligned}
	\end{equation}
	The definition of $S_\mu(a,b)$ then gives
	\begin{equation}\label{eq:strict-concavity-minimizer-proof}
		S_\mu(a,b)+o(1)
		\geq S_\mu(a,b)\left(
		\mathcal D_{\mathcal C}(w_*)^{1/p}
		+\mathcal D_{\mathcal C}(z_n)^{1/p}\right)+o(1).
	\end{equation}
	Set $a_0=\mathcal D_{\mathcal C}(w_*)>0$ and, after
	passing to a subsequence, let
	$a_1=\lim_{n\to\infty}\mathcal D_{\mathcal C}(z_n)\geq0$.
	By \eqref{eq:minimizer-splitting}, $a_0+a_1=1$.
	If $a_1>0$, the strict concavity of $t\mapsto t^{1/p}$
	yields $a_0^{1/p}+a_1^{1/p}>1$, contradicting
	\eqref{eq:strict-concavity-minimizer-proof}.
	Thus $a_1=0$ and $a_0=1$, so
	\[
	\mathcal D_{\mathcal C}(w_*)=1,\qquad
	\mathcal D_{\mathcal C}(z_n)\to0.
	\]
	
	By the definition of $S_\mu(a,b)$,
	$\|w_*\|_\kappa^2\geq S_\mu(a,b)$, whereas weak lower
	semicontinuity and \eqref{eq:minimizing-normalization}
	give the reverse inequality. Hence
	$\|w_*\|_\kappa^2=S_\mu(a,b)$, and
	\eqref{eq:minimizer-splitting} implies
	$z_n\to0$ in $H^1(\mathcal C)$.
	The inverse Emden--Fowler transform of $w_*$
	therefore attains $S_\mu(a,b)$ in
	$D_a^{1,2}(\mathbb R^N)$.
\end{proof}

\begin{proposition}[Euler equation and positivity]\label{prop:extremal-equation}
	An extremal can be chosen nonnegative.  After multiplication by a positive constant, it satisfies \eqref{eq:pure-critical-equation-intro} and \eqref{eq:critical-normalization-intro}.  It is positive almost everywhere. If $p>2$, its continuous representative is strictly positive on $\mathbb R^N\setminus\{0\}$.
\end{proposition}

\begin{proof}
	Since replacing $u$ by $|u|$ leaves $A(u)$ and
	$\mathcal D(u)$ unchanged, we may choose a nonnegative
	minimizer $W$ with $\mathcal D(W)=1$. By the Lagrange multiplier rule and testing the resulting equation with \(W\), we have
	\[
	-\operatorname{div}(|x|^{-2a}\nabla W)
	=S_\mu(a,b)|x|^{-bp}
	\left(\int_{\mathbb R^N}
	\frac{W(y)^p}{|x-y|^\mu|y|^{bp}}\,dy\right)
	W^{p-1}.
	\]
	Set $U=S_\mu(a,b)^{1/(2p-2)}W$. By homogeneity,
	$U$ satisfies \eqref{eq:pure-critical-equation-intro}
	and \eqref{eq:critical-normalization-intro}.
	
	Let $w$ be the Emden--Fowler transform of $U$.
	Since $w\geq0$, $w\not\equiv0$, and
	$L_\kappa w\geq0$ weakly in $\mathcal C$, the weak
	Harnack inequality gives a positive essential lower
	bound on a smaller ball where $w$ is nontrivial.
	Applying the same estimate along a chain of overlapping
	balls, we propagate this positivity throughout
	$\mathcal C$. Thus $w>0$ almost everywhere.
	This argument applies for every $p>1$ and does not
	require boundedness of $w^{p-2}$.
	If $p>2$, Lemma \ref{lem:cylinder-regularity} provides
	a continuous representative of $w$. The local positive
	lower bounds obtained above then imply that this
	representative is strictly positive throughout
	$\mathcal C$.
\end{proof}

In cylinder variables, the normalized extremal solves
\begin{equation}\label{eq:cylinder-critical-equation}
	-\partial_{tt}w-\Delta_{\omega}w+\frac{\kappa^2}{4}w
	=\left(\int_{\mathcal C}\mathcal K_{\mu}(t-s,\omega,\omega')|w(s,\omega')|^p\,ds\,d\omega'\right)w^{p-1}.
\end{equation}

\begin{lemma}[Regularity on the cylinder]\label{lem:cylinder-regularity}
	Let $p>2$ and let $0\leq w\in H^1(\mathcal C)$ solve \eqref{eq:cylinder-critical-equation}.  Then
	\begin{equation}\label{eq:w-bounded}
		w\in L^{\infty}(\mathcal C)\cap C^{1,\sigma}_{\mathrm{loc}}(\mathcal C)
	\end{equation}
	for some $\sigma\in(0,1)$, and
	\begin{equation}\label{eq:w-uniform-zero}
		\lim_{|t|\to\infty}\sup_{\omega\in\mathbb S^{N-1}}w(t,\omega)=0.
	\end{equation}
\end{lemma}

\begin{proof}
	Set $\ell=2^\sharp$ and $V=\mathcal K_\mu*w^p$.
	The cylindrical HLS inequality gives
	$\|V\|_{L^{2N/\mu}(\mathcal C)}
	\leq C\|w\|_{L^\ell(\mathcal C)}^p$.
	Thus the equation $L_\kappa w=Qw$, with
	$Q=Vw^{p-2}$, satisfies
	\begin{equation}\label{eq:Q-integrability-regularity-proof}
		\|Q\|_{L^{N/(2d)}(\mathcal C)}
		\leq C\|w\|_{L^\ell(\mathcal C)}^{2p-2},
		\qquad
		\frac{2d}{N}
		=\frac{\mu}{2N}+\frac{p-2}{\ell}
		<\frac2N.
	\end{equation}
	In particular, $N/(2d)>N/2$. Applying
	\eqref{eq:uniform-local-boundedness}, we obtain
	\[
	\|w\|_{L^\infty(K_1(t_0))}
	\leq C\|w\|_{L^2(K_2(t_0))}
	\qquad\text{for all }t_0\in\mathbb R,
	\]
	where $C$ is independent of $t_0$. Since
	$w\in L^2(\mathcal C)$, the right-hand side is
	uniformly bounded and tends to zero as
	$|t_0|\to\infty$. Consequently,
	$w\in L^\infty(\mathcal C)$ and
	\eqref{eq:w-uniform-zero} holds.
	
	It remains to establish the regularity of $w$.
	Since $0<\mu<N$, the kernel $\mathcal K_\mu$
	is uniformly integrable in either variable. Indeed,
	\[
	\sup_{z\in\mathcal C}
	\int_{\mathcal C}\mathcal K_\mu(z,z')\,dz'
	\leq C\left(
	\int_0^1 r^{N-1-\mu}\,dr
	+1+\int_1^\infty e^{-\mu r/2}\,dr
	\right)<\infty.
	\]
	The boundedness of $w$ therefore implies
	$V\in L^\infty(\mathcal C)$. Hence the right-hand
	side $Vw^{p-1}$ of the cylinder equation is locally
	bounded. Standard interior elliptic estimates give
	$w\in W^{2,r}_{\mathrm{loc}}(\mathcal C)$ for every
	$1\leq r<\infty$. Choosing $r>N$ and using the
	Sobolev embedding, we conclude that
	$w\in C^{1,\sigma}_{\mathrm{loc}}(\mathcal C)$
	for some $\sigma\in(0,1)$.
\end{proof}

\begin{proposition}[Decay and asymptotic behavior on the cylinder]\label{prop:axial-decay}
	Let $w$ be a positive solution in Lemma \ref{lem:cylinder-regularity}.  There exist $c,C>0$ such that
	\begin{equation}\label{eq:axial-decay}
		c e^{-\kappa|t|/2}
		\leq w(t,\omega)\leq C e^{-\kappa|t|/2},
		\qquad (t,\omega)\in\mathcal C.
	\end{equation}
Moreover,
\begin{equation}\label{eq:axial-gradient-decay}
|\partial_tw(t,\omega)|+|\nabla_\omega w(t,\omega)|
\leq Ce^{-\kappa|t|/2}.
\end{equation}

	Put $F=(\mathcal K_\mu*w^p)w^{p-1}$. The numbers
	\begin{equation}\label{eq:endpoint-coefficients}
		C_\pm=\frac{1}{\kappa|\mathbb S^{N-1}|}
		\int_{\mathcal C}e^{\pm\kappa s/2}F(s,\omega')\,ds\,d\omega'
	\end{equation}
	are finite and positive, and
	\begin{equation}\label{eq:axial-endpoint-limits}
		\begin{aligned}
			\sup_\omega\bigl|e^{\kappa t/2}w(t,\omega)-C_+\bigr|&\longrightarrow0
			&& (t\to+\infty),\\
			\sup_\omega\bigl|e^{-\kappa t/2}w(t,\omega)-C_-\bigr|&\longrightarrow0
			&& (t\to-\infty).
		\end{aligned}
	\end{equation}
\end{proposition}

\begin{proof}
	Set $\nu=\kappa/2$, $V=\mathcal K_\mu*w^p$ and
	$F=Vw^{p-1}$. Then $L_\kappa w=F$, where
	$L_\kappa=-\partial_{tt}-\Delta_{\mathbb S^{N-1}}+\nu^2$.
	
	We first recall the estimates for the Green kernel of
	$L_\kappa$. It admits the representation
	\[
	G_\kappa(t,\omega;s,\omega')
	=\int_0^\infty e^{-\nu^2r}(4\pi r)^{-1/2}
	e^{-(t-s)^2/(4r)}
	H_{\mathbb S^{N-1}}(r,\omega,\omega')\,dr.
	\]
	The Gaussian heat kernel bound of Li and Yau
	\cite[Corollary~3.1]{LiYau1986}, together with the
	boundedness of $H_{\mathbb S^{N-1}}$ for $r\geq1$,
	gives
	\[
	G_\kappa(z,z')\leq C\bigl(1+d_{\mathcal C}(z,z')^{2-N}\bigr),
	\qquad z\neq z'.
	\]

	    We next derive the spherical harmonic expansion of
	the Green kernel. Let $\{Y_{\ell,m}\}_{m=1}^{d_\ell}$
	be a real orthonormal basis of spherical harmonics
	of degree $\ell$ on $\mathbb S^{N-1}$. Set
	$\nu_\ell=\sqrt{\ell(\ell+N-2)+\nu^2}$.
	Expanding $G_\kappa$ in spherical harmonics, we find
	that each coefficient solves
	\[
	\left(-\frac{d^2}{d\tau^2}+\nu_\ell^2\right)
	g_\ell(\tau)=\delta_0,
	\qquad \tau=t-s.
	\]
	The decaying solution is
	$g_\ell(\tau)=e^{-\nu_\ell|\tau|}/(2\nu_\ell)$.
	Consequently,
	\begin{equation}\label{eq:green-spherical-expansion}
		G_\kappa(t,\omega;s,\omega')
		=\sum_{\ell=0}^\infty\sum_{m=1}^{d_\ell}
		\frac{e^{-\nu_\ell|t-s|}}{2\nu_\ell}
		Y_{\ell,m}(\omega)Y_{\ell,m}(\omega').
	\end{equation}
	The expansion holds in the distributional sense and
	converges absolutely for $t\neq s$.
	The degree-zero mode is
	$G_0(t-s)=e^{-\nu|t-s|}/(2\nu|\mathbb S^{N-1}|)$.
	Writing $G_\kappa^\perp=G_\kappa-G_0$, the addition
	formula for spherical harmonics gives
	\[
	\begin{aligned}
		|G_\kappa^\perp|
		\leq C\sum_{\ell=1}^\infty
		\frac{(1+\ell)^{N-2}}{\nu_\ell}
		e^{-\nu_\ell|t-s|}
		\leq Ce^{-\nu_1|t-s|}
		\sum_{\ell=1}^\infty
		\frac{(1+\ell)^{N-2}}{\nu_\ell}
		e^{-(\nu_\ell-\nu_1)}
		\leq Ce^{-\nu_1|t-s|}
	\end{aligned}
	\]
	for $|t-s|\geq1$, where
	$\nu_1=\sqrt{N-1+\nu^2}>\nu$.
	The series converges since $\nu_\ell$ grows linearly
	in $\ell$. Combining this estimate with the local
	Green kernel bound, we obtain
	\[
	|G_\kappa^\perp(z,z')|
	\leq
	\begin{cases}
		C(1+d_{\mathcal C}(z,z')^{2-N}),
		& |t-s|<1,\\
		Ce^{-\nu_1|t-s|},&|t-s|\geq1.
	\end{cases}
	\]
	
	We next establish a preliminary decay estimate for $w$.
	By Lemma~\ref{lem:cylinder-regularity}, $V$ is bounded
	and $w\to0$ uniformly in $\omega$ as $|t|\to\infty$.
	Since $p>2$, for every $\varepsilon>0$ we have
	$F\leq\varepsilon w$ whenever $|t|$ is sufficiently large.
	Fix $0<\sigma<\nu$ and choose
	$0<\varepsilon<\nu^2-\sigma^2$. For some sufficiently
	large $T$, it follows that
	$(L_\kappa-\varepsilon)w\leq0$ on
	$\mathcal C_T=(T,\infty)\times\mathbb S^{N-1}$.
	The function $\Phi(t)=Ce^{-\sigma(t-T)}$, with $C$
	sufficiently large, satisfies $w\leq\Phi$ at $t=T$
	and
	\[
	(L_\kappa-\varepsilon)\Phi
	=(\nu^2-\varepsilon-\sigma^2)\Phi\geq0.
	\]
	Since $w\leq\Phi$ at $t=T$, the function
	$h=(w-\Phi)_+$ has zero trace on the boundary of
	$\mathcal C_T$ and belongs to $H_0^1(\mathcal C_T)$.
	Testing the inequality for $w-\Phi$ with
	$h=(w-\Phi)_+$, by a standard cutoff argument, we obtain
	\[
	\int_{\mathcal C_T}
	\left(|\nabla_{\mathcal C}(w-\Phi)_+|^2
	+(\nu^2-\varepsilon)(w-\Phi)_+^2\right)
	\,dt\,d\omega\leq0.
	\]
	Hence $w\leq\Phi$. The same comparison on the left
	half-cylinder, together with boundedness on compact
	slabs, yields
	\[
	w(t,\omega)\leq C_\sigma e^{-\sigma|t|}
	\qquad\text{for every }0<\sigma<\nu.
	\]
	
	Choose $\sigma$ such that $\nu/(p-1)<\sigma<\nu$ and
	set $b_1=(p-1)\sigma>\nu$. Since $V$ is bounded,
	\begin{equation}\label{eq:F-weighted-integrable}
		0\leq F(t,\omega)\leq Ce^{-b_1|t|},
		\qquad
		\int_{\mathcal C}e^{\nu|s|}
		F(s,\omega')\,ds\,d\omega'<\infty.
	\end{equation}
	In particular, $F\in L^2(\mathcal C)$. Coercivity of
	$L_\kappa$ and uniqueness for the linear equation
	$L_\kappa w=F$ therefore give the Green representation
	\[
	w(t,\omega)=\int_{\mathcal C}
	G_\kappa(t,\omega;s,\omega')
	F(s,\omega')\,ds\,d\omega'.
	\]
	
    We now use this decomposition to derive the precise
    asymptotics of $w$.
	Choose $\nu<\sigma_1<\min\{\nu_1,b_1\}$. On
	$|t-s|<1$, the weighted estimate
	\eqref{eq:F-weighted-integrable} gives
	$F(s,\omega')\leq Ce^{-b_1|t|}$. Moreover, the local
	Green kernel singularity is integrable uniformly in
	$(t,\omega)$. Splitting the integral into
	$|t-s|<1$ and $|t-s|\geq1$, we obtain
	\[
		\left|\int_{\mathcal C}
		G_\kappa^\perp(t,\omega;s,\omega')
		F(s,\omega')\,ds\,d\omega'\right|
		\leq
		Ce^{-b_1|t|}
		+C\int_{\mathbb R}e^{-\nu_1|t-s|}
		e^{-b_1|s|}\,ds
		\leq Ce^{-\sigma_1|t|}.
	\]
	This estimate is uniform in $\omega$.
	
	It remains to consider the constant spherical mode,
	\[
	w_0(t)=\frac{1}{2\nu|\mathbb S^{N-1}|}
	\int_{\mathcal C}e^{-\nu|t-s|}
	F(s,\omega')\,ds\,d\omega'.
	\]
	As $t\to+\infty$, the factor
	$e^{\nu t-\nu|t-s|}$ converges to $e^{\nu s}$
	and is bounded by $e^{\nu|s|}$. By
	\eqref{eq:F-weighted-integrable} and dominated
	convergence,
	\[
	\lim_{t\to+\infty}e^{\nu t}w_0(t)
	=\frac{1}{2\nu|\mathbb S^{N-1}|}
	\int_{\mathcal C}e^{\nu s}F(s,\omega')\,ds\,d\omega'
	=C_+.
	\]
	Similarly,
	\[
	\lim_{t\to-\infty}e^{-\nu t}w_0(t)
	=\frac{1}{2\nu|\mathbb S^{N-1}|}
	\int_{\mathcal C}e^{-\nu s}F(s,\omega')\,ds\,d\omega'
	=C_-.
	\]
	Since $F\geq0$ and $F\not\equiv0$, both constants
	are finite and strictly positive. The uniform estimate
	for $G_\kappa^\perp$ and $\sigma_1>\nu$ now yield
	\eqref{eq:axial-endpoint-limits}, with
	$2\nu=\kappa$ giving precisely
	\eqref{eq:endpoint-coefficients}.
	
	The endpoint limits imply the two-sided estimate
	\eqref{eq:axial-decay} for sufficiently large $|t|$.
	Positivity and continuity of $w$ extend this estimate
	to the whole cylinder.
	
	To estimate the derivatives of $w$, we fix $r>N$ and
	set $K_j(\tau)=(\tau-j,\tau+j)\times\mathbb S^{N-1}$
	for $j=1,2$. Applying the interior elliptic estimate
	and using the decay of $w$ and $F$, we obtain
	\[
	\|w\|_{W^{2,r}(K_1(\tau))}
	\leq C\left(
	\|w\|_{L^r(K_2(\tau))}
	+\|F\|_{L^r(K_2(\tau))}\right)
	\leq Ce^{-\nu|\tau|}.
	\]
	The constant is independent of $\tau$ by axial
	translation invariance. We then use the Sobolev embedding to deduce
	\eqref{eq:axial-gradient-decay}.
\end{proof}

\begin{proof}[Proof of Theorem \ref{thm:critical-extremal}]
	By Propositions \ref{prop:attainment} and
	\ref{prop:extremal-equation}, the quotient
	$S_\mu(a,b)$ is attained by a positive extremal $U$
	satisfying \eqref{eq:pure-critical-equation-intro}
	and \eqref{eq:critical-normalization-intro}.
	Assume $p>2$. Writing $x=r\omega$ with $r=|x|$, we have
	$U(r\omega)=r^{-\kappa/2}w(-\log r,\omega)$.
	Applying Proposition \ref{prop:axial-decay} to this
	representation, we obtain
	\[
	\lim_{r\downarrow0}\sup_{\omega\in\mathbb S^{N-1}}
	|U(r\omega)-C_+|=0,
	\qquad
	\lim_{r\to\infty}\sup_{\omega\in\mathbb S^{N-1}}
	|r^\kappa U(r\omega)-C_-|=0.
	\]
	The two-sided bounds \eqref{eq:axial-decay} also give
	$c\leq U(x)\leq C$ for $0<|x|\leq1$ and
	$c|x|^{-\kappa}\leq U(x)\leq C|x|^{-\kappa}$
	for $|x|\geq1$. These estimates establish
	\eqref{eq:extremal-decay-intro} and
	\eqref{eq:extremal-endpoint-limits-intro}.
\end{proof}

\begin{corollary}[Finite mass of a critical extremal]\label{cor:finite-mass-extremal}
	Under \eqref{eq:decay-range}, the extremal in Theorem \ref{thm:critical-extremal} belongs to $X$ if and only if $\kappa>2$.
\end{corollary}
\begin{proof}
	Since $U$ is bounded near the origin and
	$N-2a=\kappa+2>0$, we have
	\[
	\int_{|x|<1}|x|^{-2a}U^2\,dx
	\leq C\int_0^1r^{N-1-2a}\,dr<\infty.
	\]
	At infinity, the two-sided decay estimate
	\eqref{eq:extremal-decay-intro} gives
	\[
	\int_{|x|>1}|x|^{-2a}U^2\,dx
	\asymp\int_1^\infty r^{1-\kappa}\,dr,
	\]
	which is finite if and only if $\kappa>2$.
	Since $U\in D_a^{1,2}(\mathbb R^N)$, we conclude
	that $U\in X$ precisely when $\kappa>2$.
\end{proof}

We next estimate truncations of the extremal.  For $R>e$, choose $\chi_R\in C_c^{\infty}(\mathbb R)$ such that $0\leq\chi_R\leq1$,
\begin{equation}\label{eq:cylinder-cutoff-properties}
	\chi_R(t)=1\quad\text{if }|t|\leq\log R,
	\qquad
	\chi_R(t)=0\quad\text{if }|t|\geq\log(2R),
\end{equation}
and $|\chi_R'|+|\chi_R''|\leq C$, with $C$ independent of $R$.  Put
\[
w_R(t,\omega)=\chi_R(t)w(t,\omega)
\]
and let $U_R$ be its inverse Emden--Fowler transform.

\begin{corollary}[Truncated extremals]\label{cor:truncated-extremals}
	Assume \eqref{eq:decay-range}, equivalently $p>2$. The truncated functions $U_R\in X$ are nonnegative and satisfy
	\[
	\operatorname{supp}U_R\subset\{(2R)^{-1}\leq|x|\leq2R\},
	\]
	\[
	A(U_R)=A(U)+O(R^{-\kappa}),
	\qquad
	\mathcal D(U_R)=\mathcal D(U)+O(R^{-\kappa}),
	\]
	and
	\begin{equation}\label{eq:truncated-mass}
		M(U_R)\leq C
		\begin{cases}
			1,&\kappa>2,\\
			1+\log R,&\kappa=2,\\
			R^{2-\kappa},&0<\kappa<2.
		\end{cases}
	\end{equation}
	For every $q\in(2,2^{\sharp})$,
	\begin{equation}\label{eq:truncated-B-lower}
		B_q(U_R)\geq c_q>0
	\end{equation}
	for all large $R$.
\end{corollary}

\begin{proof}
	Combining \eqref{eq:axial-decay} and
	\eqref{eq:axial-gradient-decay}, we obtain
	\[
	\int_{|t|>\log R}\int_{\mathbb S^{N-1}}
	\left(|\partial_tw|^2+|\nabla_\omega w|^2+w^2\right)
	\,d\omega\,dt=O(R^{-\kappa}).
	\]
	Since $w_R=\chi_Rw$ and $\chi_R$ depends only on $t$,
	we have $\partial_tw_R=\chi_R\partial_tw+\chi_R'w$
	and $\nabla_\omega w_R=\chi_R\nabla_\omega w$.
	Expanding the quadratic form gives
	\[
	\begin{aligned}
		A(U_R)-A(U)
		={}&\int_{\mathcal C}(\chi_R^2-1)
		\left(|\partial_tw|^2+|\nabla_\omega w|^2
		+\frac{\kappa^2}{4}w^2\right)\,dt\,d\omega\\
		&+2\int_{\mathcal C}\chi_R\chi_R'w\partial_tw
		\,dt\,d\omega
		+\int_{\mathcal C}|\chi_R'|^2w^2\,dt\,d\omega.
	\end{aligned}
	\]
	All three terms are supported in $|t|>\log R$.
	The decay estimates and the uniform bounds on $\chi_R$
	and $\chi_R'$ therefore give
	$A(U_R)=A(U)+O(R^{-\kappa})$.
	
	For the nonlocal term, set $f=|w|^p$ and
	$f_R=|w_R|^p$. Since $pr_\mu=2^\sharp$, the decay
	estimate yields
	\[
	\|f-f_R\|_{r_\mu}^{r_\mu}
	\leq C\int_{|t|>\log R}
	e^{-\kappa pr_\mu|t|/2}\,dt
	\leq CR^{-\kappa pr_\mu/2}.
	\]
	Hence $\|f-f_R\|_{r_\mu}\leq CR^{-\kappa p/2}$.
	By the cylindrical HLS inequality,
	\[
		|\mathcal D(U_R)-\mathcal D(U)|
		\leq C\|f-f_R\|_{r_\mu}
		\bigl(\|f\|_{r_\mu}+\|f_R\|_{r_\mu}\bigr)
		=O(R^{-\kappa}),
	\]
	where we use $p>2$ in the last estimate.
	
	We next estimate the mass. By
	\eqref{eq:cylinder-mass-local} and
	\eqref{eq:axial-decay}, we have
	\[
	M(U_R)\leq C\left(
	1+\int_{-\log(2R)}^0e^{(\kappa-2)t}\,dt\right).
	\]
    According to the value of $\kappa$, we obtain
	\[
	M(U_R)\leq C
	\begin{cases}
		1,&\kappa>2,\\
		1+\log R,&\kappa=2,\\
		R^{2-\kappa},&0<\kappa<2,
	\end{cases}
	\]
	which proves \eqref{eq:truncated-mass}.
	
	To show that $B_q(U_R)$ stays bounded away from zero,
	choose a fixed annulus
	$\mathcal A=\{x:r_0<|x|<r_1\}$ with
	$0<r_0<r_1<\infty$. For all sufficiently large $R$,
	the cutoff is identically one on $\mathcal A$. Hence
	\[
	B_q(U_R)\geq\int_{\mathcal A}|x|^{-bq}U^q\,dx=:c_q>0,
	\]
	where $c_q$ is independent of $R$.
\end{proof}

\medskip

We next describe the loss of compactness of Palais--Smale sequences.

\subsection{Pure and constrained global compactness}\label{sec:global}

We first prove the profile decomposition for the limiting critical problem.  In cylinder variables, let
\[
\mathcal J(w)=\frac12\|w\|_{\kappa}^2-\frac1{2p}\mathcal D_{\mathcal C}(w).
\]

\begin{lemma}[One-profile extraction]\label{lem:one-profile}
	Let $z_n\rightharpoonup0$ in $H^1(\mathcal C)$ and assume
	\begin{equation}\label{eq:remainder-PS}
		\mathcal J'(z_n)\to0,
		\qquad
		\liminf_{n\to\infty}\|z_n\|_{\kappa}>0.
	\end{equation}
	Then there exist translations $s_n\in\mathbb R$, $|s_n|\to\infty$, and a nontrivial critical point $W$ of $\mathcal J$ such that
	\[
	z_n(t+s_n,\omega)\rightharpoonup W(t,\omega)
	\quad\text{in }H^1(\mathcal C).
	\]
	If $r_n=z_n-W(\cdot-s_n,\cdot)$, then
	\begin{equation}\label{eq:one-profile-remainder-PS}
		\mathcal J'(r_n)\to0\quad\text{in }H^{-1}(\mathcal C),
	\end{equation}
	and
	\begin{equation}\label{eq:one-profile-splitting}
		\begin{aligned}
			\|z_n\|_{\kappa}^2
			&=\|W\|_{\kappa}^2+\|r_n\|_{\kappa}^2+o(1),\\
			\mathcal J(z_n)&=\mathcal J(W)+\mathcal J(r_n)+o(1).
		\end{aligned}
	\end{equation}
\end{lemma}

\begin{proof}
	Testing $\mathcal J'(z_n)$ with $z_n$, we obtain
	\[
	\|z_n\|_\kappa^2-\mathcal D_{\mathcal C}(z_n)=o(1).
	\]
	Since $\|z_n\|_\kappa$ is bounded away from zero,
	we have $\mathcal D_{\mathcal C}(z_n)\geq c>0$ for
	large $n$. The cylindrical HLS inequality
	\eqref{eq:D-cylinder-bound} then shows that
	$\|z_n\|_{2^\sharp}$ is bounded away from zero.
	Hence the vanishing condition in Lemma
	\ref{lem:lions-cylinder} cannot hold. We can therefore
	choose translations $s_n\in\mathbb R$ and a subsequence
	such that
	\[
	\int_{(-1,1)\times\mathbb S^{N-1}}
	|z_n(t+s_n,\omega)|^2\,dt\,d\omega
	\geq c_0>0.
	\]
	After passing to a further subsequence, we have
	$z_n(\cdot+s_n,\cdot)\rightharpoonup W$ in
	$H^1(\mathcal C)$. The local $L^2$ lower bound and
	local compactness imply that $W\not\equiv0$. We also have $|s_n|\to\infty$. Otherwise, a subsequence
	would satisfy $s_n\to s$ for some $s\in\mathbb R$.
	Since $z_n\rightharpoonup0$ in $H^1(\mathcal C)$,
	the translated sequence would then converge weakly
	to zero, contradicting $W\not\equiv0$.
	
	Set $\widetilde z_n=z_n(\cdot+s_n,\cdot)$.
	Translation invariance gives
	$\mathcal J'(\widetilde z_n)\to0$.
	Testing against compactly supported functions and using
	local compactness together with Lemma
	\ref{lem:nonlocal-BL}, we obtain $\mathcal J'(W)=0$.
	Set $r_n=z_n-W(\cdot-s_n,\cdot)$. Weak convergence
	and Lemma \ref{lem:nonlocal-BL} give the splittings
	\eqref{eq:one-profile-splitting}. Applying the derivative
	splitting \eqref{eq:D-BL-derivative} to $\widetilde z_n$
	and using $\mathcal J'(W)=0$, we obtain
	$\mathcal J'(\widetilde z_n-W)\to0$ in
	$H^{-1}(\mathcal C)$. Since axial translations preserve the $H^{-1}(\mathcal C)$
	norm, we conclude that
	$\mathcal J'(r_n)\to0$ in $H^{-1}(\mathcal C)$.
\end{proof}
\begin{proof}[Proof of Theorem \ref{thm:global-compactness-intro}(i)]
	Let \(w_n\) be the cylinder transforms of \(u_n\). Up to a subsequence, we have
	\[
	w_n\rightharpoonup w^0\quad\text{in }H^1(\mathcal C).
	\]
	By local compactness and Lemma \ref{lem:nonlocal-BL}, we have \(\mathcal J'(w^0)=0\). Set \(z_n^1=w_n-w^0\). The derivative splitting implies that \(\mathcal J'(z_n^1)\to0\). If \(z_n^1\to0\) in \(H^1(\mathcal C)\), the proof is complete. Otherwise, we apply Lemma \ref{lem:one-profile} to extract a nontrivial profile. Repeating this procedure, we obtain profiles \(W^j\) and translations \(s_n^j\).
	
	We next show that the translations satisfy
	\[
	|s_n^i-s_n^j|\to\infty\qquad(i\neq j).
	\]
	Let \(r_n^{j-1}\) denote the remainder before the \(j\)-th extraction. By construction,
	\[
	r_n^{j-1}(\cdot+s_n^i,\cdot)\rightharpoonup0
	\quad(i<j),
	\qquad
	r_n^{j-1}(\cdot+s_n^j,\cdot)\rightharpoonup W^j\neq0.
	\]
	Suppose that \(|s_n^j-s_n^i|\) fails to tend to infinity. Passing to a subsequence, we may assume that $s_n^j-s_n^i\to\tau\in\mathbb R$. For every \(\varphi\in C_c^\infty(\mathcal C)\), we have
	\[
	\langle r_n^{j-1}(\cdot+s_n^j),\varphi\rangle_\kappa
	=\langle r_n^{j-1}(\cdot+s_n^i),
	\varphi(\cdot-(s_n^j-s_n^i))\rangle_\kappa
	\longrightarrow0,
	\]
	since the translated test functions converge strongly in \(H^1(\mathcal C)\). This contradicts \(W^j\neq0\). Hence
	\begin{equation}\label{eq:translation-orthogonality-proof}
		|s_n^i-s_n^j|\to\infty\qquad(i\neq j).
	\end{equation}
	
	For any nontrivial critical point \(W\) of \(\mathcal J\), testing \(\mathcal J'(W)=0\) with \(W\) and using the definition of \(S_\mu(a,b)\), we have
	\[
	\|W\|_\kappa^2=\mathcal D_{\mathcal C}(W)>0,
	\qquad
	\|W\|_\kappa^2\geq S_\mu(a,b)\|W\|_\kappa^{2/p}.
	\]
	Thus
	\[
	\|W\|_\kappa^2\geq S_\mu(a,b)^{p/(p-1)},
	\qquad
	\mathcal J(W)\geq T_{\mu,a,b}.
	\]
	After extracting \(j\) profiles, we have
	\[
	jS_\mu(a,b)^{p/(p-1)}
	\leq\sum_{i=1}^j\|W^i\|_\kappa^2
	\leq\limsup_n\|w_n\|_\kappa^2-\|w^0\|_\kappa^2.
	\]
	Since \((w_n)\) is bounded in \(H^1(\mathcal C)\), the extraction process terminates after finitely many steps.
	
	Let \(r_n\) be the final remainder. If \(r_n\not\to0\) in \(H^1(\mathcal C)\), then \(\mathcal J'(r_n)[r_n]=o(1)\), \eqref{eq:D-cylinder-bound}, and Lemma \ref{lem:lions-cylinder} imply the existence of another nontrivial profile, a contradiction. Thus \(r_n\to0\) in \(H^1(\mathcal C)\).
	
	Finally, setting \(\varepsilon_n^j=e^{-s_n^j}\), we obtain \eqref{eq:scale-orthogonality-intro} from \eqref{eq:translation-orthogonality-proof}. Returning to the Euclidean variables, the standard Hilbert space decomposition and nonlocal Brézis--Lieb splittings, we obtain \eqref{eq:profile-decomposition-intro}--\eqref{eq:profile-splitting-intro}.
\end{proof}

We next prove the compactness theorem for constrained sequences.
\begin{proof}[Proof of Theorem \ref{thm:global-compactness-intro}(ii)]
	Set $\nu=\kappa/2$.
	
	Since \(\|d(E_\beta|_{S_\rho})(u_n)\|\to0\), there exist \(\lambda_n\in\mathbb R\) such that
	\begin{equation}\label{eq:EL-PS-sequence}
		E_\beta'(u_n)-\frac{\lambda_n}{2}M'(u_n)\to0\quad\text{in }X^*.
	\end{equation}
	Testing \eqref{eq:EL-PS-sequence} with \(u_n\) and using \(P_\beta(u_n)=o(1)\), we obtain
	\begin{equation}\label{eq:lambda-n-formula}
		\lambda_n\rho^2=\beta(\delta_q-1)B_q(u_n)+o(1).
	\end{equation}
	Since \((u_n)\) is bounded in \(X\), there exists a subsequence still denoted by \((u_n)\) such that $u_n\rightharpoonup u^0$ in $X$. By Proposition \ref{prop:compact-embedding}, we have
	\begin{equation}\label{eq:Bq-constrained-compact-proof}
		B_q(u_n)\to B_q(u^0),\qquad
		\lambda_n\to\lambda\leq0.
	\end{equation}
	By local compactness and Lemma \ref{lem:nonlocal-BL}, passing to the limit in \eqref{eq:EL-PS-sequence}, we have
	\[
	\langle E_\beta'(u^0),\varphi\rangle
	-\lambda\int_{\mathbb R^N}|x|^{-2a}u^0\varphi\,dx=0,
	\qquad
	\varphi\in C_c^\infty(\mathbb R^N\setminus\{0\}).
	\]
	
	To extend this identity to \(X\), let \(\chi_R\) be a smooth cutoff satisfying
	\(0\leq\chi_R\leq1, \chi_R=1 \text{ on }\{R^{-1}\leq|x|\leq R\}, \chi_R=0 \text{ on } \{|x|\leq(2R)^{-1}\}\cup\{|x|\geq2R\}, |\nabla\chi_R|\leq C|x|^{-1}.\)
	For \(\varphi\in X\), set \(\varphi_R=\chi_R\varphi\). By the weighted Hardy inequality, we have
	\[
	A(\varphi-\varphi_R)
	\leq C\int_{\{|x|<R^{-1}\}\cup\{|x|>R\}}
	\left(
	|x|^{-2a}|\nabla\varphi|^2
	+|x|^{-2a-2}\varphi^2
	\right)\,dx
	\longrightarrow0.
	\]
	Moreover, by dominated convergence, we have
	\[
	M(\varphi-\varphi_R)\longrightarrow0.
	\]
	Thus \(\varphi_R\to\varphi\) in \(X\). Since \(\varphi_R\) has compact support away from the origin, it can be approximated in \(X\) by smooth compactly supported functions. Hence \(C_c^\infty(\mathbb R^N\setminus\{0\})\) is dense in \(X\). By density and continuity, \(u^0\) satisfies the Euler equation with Lagrange multiplier \(\lambda\) in \(X^*\).
	We first prove that a nonzero weak limit has mass $\rho^2$. If $u^0\neq0$, Proposition \ref{prop:pohozaev} and \eqref{eq:lambda-n-formula} give
	\[
	\lambda M(u^0)=\beta(\delta_q-1)B_q(u^0)
	=\lambda\rho^2<0.
	\]
	Hence
	\begin{equation}\label{eq:mass-recovery}
		\lambda<0,\qquad M(u^0)=\rho^2,\qquad M(u_n-u^0)\to0.
	\end{equation}
	Indeed, weak convergence and equality of the limiting mass give $M(u_n-u^0)=M(u_n)-M(u^0)+o(1)=o(1)$. If $u^0=0$, then \eqref{eq:lambda-n-formula} gives $\lambda=0$. In both cases,
	\begin{equation}\label{eq:weak-limit-pohozaev-proof}
		P_\beta(u^0)=0.
	\end{equation}
	
	Let $w_n,w^0$ be the cylinder transforms and put $z_n^0=w_n-w^0$. The profiles are extracted in $H^1(\mathcal C)$. If $z_n^0\to0$ in $L^{2^\sharp}$, there is no profile to extract. Otherwise Lemma \ref{lem:lions-cylinder} gives shifts $s_n$ and a nonzero weak limit of $z_n^0(\cdot+s_n)$. Local compactness before translation forces $|s_n|\to\infty$.
	
	We show that the shifts cannot tend to $-\infty$. If $s_n\to-\infty$ and $w_n(\cdot+s_n)\rightharpoonup W\neq0$, choose a fixed slab $K_L$ with $\int_{K_L}|W|^2\,dt\,d\omega>0$. Local strong convergence yields
	\begin{equation}\label{eq:no-minus-infinity-profile}
		\rho^2=M(u_n)
		\geq e^{-2s_n-2L}\int_{K_L}|w_n(t+s_n,\omega)|^2\,dt\,d\omega
		\longrightarrow\infty,
	\end{equation}
	a contradiction. Since $w^0(\cdot+s_n)\to0$ locally as $|s_n|\to\infty$, the same conclusion applies to the first remainder. Thus $s_n\to+\infty$.
	
	Let \(\varphi\in C_c^\infty(\mathcal C)\) and set
	\[
	\widetilde w_n(t,\omega)=w_n(t+s_n,\omega).
	\]
	Let $\psi_n\in X$ denote the function whose cylindrical transform is $\varphi(\cdot-s_n,\cdot)$. Then
	\[
	\|\psi_n\|_X^2
	=\|\varphi\|_\kappa^2
	+e^{-2s_n}\int_{\mathcal C}e^{-2t}|\varphi|^2\,dt\,d\omega
	\leq C.
	\]
	By \eqref{eq:EL-PS-sequence} and the uniform boundedness of $\{\psi_n\}$ in $X$, the corresponding error term satisfies $\langle E_\beta'(u_n)-\frac{\lambda_n}{2}M'(u_n),\psi_n\rangle=o(1)$.
	
	Since \((\widetilde w_n)\) is bounded in \(H^1(\mathcal C)\), \((\lambda_n)\) is bounded, and \(\gamma_q>0\), we have
	\[
	\begin{aligned}
		\left|\lambda_n\int_{\mathcal C}
		e^{-2t}w_n(t,\omega)\varphi(t-s_n,\omega)\,dt\,d\omega\right|
		&\leq Ce^{-2s_n}\longrightarrow0,\\
		\left|\beta\int_{\mathcal C}
		e^{-\gamma_qt}|w_n|^{q-2}w_n
		\varphi(t-s_n,\omega)\,dt\,d\omega\right|
		&\leq Ce^{-\gamma_qs_n}\longrightarrow0.
	\end{aligned}
	\]
	By local compactness and Lemma \ref{lem:nonlocal-BL}, we have
	\[
	\int_{\mathcal C}
	(\mathcal K_\mu*|\widetilde w_n|^p)
	|\widetilde w_n|^{p-2}\widetilde w_n\varphi\,dt\,d\omega
	\longrightarrow
	\int_{\mathcal C}
	(\mathcal K_\mu*|W|^p)|W|^{p-2}W\varphi\,dt\,d\omega.
	\]
	Passing to the limit in \eqref{eq:EL-PS-sequence}, we obtain
	\[
	\int_{\mathcal C}
	\bigl(\nabla_{\mathcal C}W\cdot\nabla_{\mathcal C}\varphi
	+\nu^2W\varphi\bigr)\,dt\,d\omega
	=
	\int_{\mathcal C}
	(\mathcal K_\mu*|W|^p)|W|^{p-2}W\varphi\,dt\,d\omega.
	\]
	By density, \(\mathcal J'(W)=0\) in \(H^{-1}(\mathcal C)\).
	\begin{equation}\label{eq:profile-pure-equation-constrained-proof}
		L_\kappa W=(\mathcal K_\mu*|W|^p)|W|^{p-2}W.
	\end{equation}
	
	After $j$ extractions, set
	\[
	z_n^j=w_n-w^0-\sum_{i=1}^jW^i(\cdot-s_n^i,\cdot).
	\]
	By construction, we have
	\[
	z_n^j\rightharpoonup0,
	\qquad
	z_n^j(\cdot+s_n^i,\cdot)\rightharpoonup0
	\quad(1\leq i\leq j)
	\]
	in $H^1(\mathcal C)$.
	If $z_n^j\not\to0$ in $L^{2^\sharp}(\mathcal C)$,
	we apply Lemma \ref{lem:lions-cylinder} to find
	translations $s_n$ such that, up to a subsequence,
	\[
	z_n^j(\cdot+s_n,\cdot)\rightharpoonup W\neq0
	\quad\text{in }H^1(\mathcal C).
	\]
	The preceding weak convergence properties imply
	\[
	|s_n|\to\infty,
	\qquad
	|s_n-s_n^i|\to\infty
	\quad(1\leq i\leq j).
	\]
	Consequently,
	\[
	w_n(\cdot+s_n,\cdot)\rightharpoonup W.
	\]
	By \eqref{eq:no-minus-infinity-profile}, we have \(s_n\to+\infty\). Using the translated test functions in \eqref{eq:EL-PS-sequence} and passing to the limit as above, we obtain \eqref{eq:profile-pure-equation-constrained-proof} for \(W\). Here the original sequence \(w_n\) is used, since the remainders need not belong to the weighted mass space.
	
	Combining $\mathcal J'(W^i)=0$ with the definition
	of $S_\mu(a,b)$, we obtain
	\[
	\|W^i\|_\kappa^2=\mathcal D_{\mathcal C}(W^i)
	\geq S_\mu(a,b)^{p/(p-1)}.
	\]
	Using the norm splitting after $j$ extractions, we obtain
	\[
	\sum_{i=1}^j\|W^i\|_\kappa^2
	\leq\limsup_n A(u_n)-A(u^0).
	\]
	Thus extraction terminates after finitely many steps with
	$z_n^k\to0$ in $L^{2^\sharp}$. In particular, $\mathcal D_{\mathcal C}(z_n^k)\to0$.
	Combining the norm and nonlocal decompositions with \eqref{eq:Bq-constrained-compact-proof} and \eqref{eq:weak-limit-pohozaev-proof}, we deduce that
	\[
	o(1)=P_\beta(u_n)
	=P_\beta(u^0)+\sum_{i=1}^k
	\bigl(\|W^i\|_\kappa^2-\mathcal D_{\mathcal C}(W^i)\bigr)
	+\|z_n^k\|_\kappa^2-\mathcal D_{\mathcal C}(z_n^k)+o(1).
	\]
	Hence $z_n^k\to0$ in $H^1$. Taking $\varepsilon_n^i=e^{-s_n^i}\to0$ proves the decomposition and all norm and energy splittings.
	
	If $u^0=0$ and $k=0$, the strong convergence $z_n^k\to0$ in $H^1(\mathcal C)$ gives $A(u_n)\to0$. Since $B_q(u_n)\to0$ and $\mathcal D(u_n)\to0$, we obtain $c=0$. If $u^0=0$ and $k\geq1$, at least one nontrivial profile contributes energy at least $T_{\mu,a,b}$, so $c\geq T_{\mu,a,b}$. If $u^0\neq0$, the mass recovery argument gives $M(u^0)=\rho^2$. These alternatives complete the proof.
\end{proof}

\section{Normalized solutions}\label{sec:normalized-solutions}

Throughout this section assume \eqref{eq:assumption-ab}, \eqref{eq:decay-range}, $\rho,\beta>0$, and $2<q<2^\sharp$. We first construct the test functions used in the energy estimates.

\subsection{Test functions and truncation estimates}\label{sec:test-functions}

As the critical extremals given by Corollary \ref{cor:finite-mass-extremal} can fail to have finite mass, we employ the truncations introduced in Corollary \ref{cor:truncated-extremals}. When working on the cylinder, we retain the notation for $A, M, B_q, E_\beta$, and $S_\rho$ to denote the corresponding Emden--Fowler images. Let $U$ be the normalized extremal from Theorem \ref{thm:critical-extremal}, and let $U_R$ be the truncation defined in Corollary \ref{cor:truncated-extremals}. Replacing $U_R$ by $c_R U_R$, where
\[
c_R=\left(\frac{A(U_R)}{\mathcal D(U_R)}\right)^{1/(2p-2)}=1+O(R^{-\kappa}),
\]
we assume
\begin{equation}\label{eq:UR-A-D}
	A(U_R)=\mathcal D(U_R).
\end{equation}
Then
\[
\frac{p-1}{2p}A(U_R)=T_{\mu,a,b}+O(R^{-\kappa}).
\]
For $\varepsilon>0$, define
\[
U_{R,\varepsilon}(x)=\varepsilon^{-\kappa/2}U_R(x/\varepsilon).
\]
A change of variables gives
\[
\begin{aligned}
	A(U_{R,\varepsilon})&=A(U_R),\qquad
	\mathcal D(U_{R,\varepsilon})=\mathcal D(U_R),\\
	M(U_{R,\varepsilon})&=\varepsilon^2M(U_R),\qquad
	B_q(U_{R,\varepsilon})=\varepsilon^{\gamma_q}B_q(U_R).
\end{aligned}
\]

\begin{lemma}[Choice of truncation and concentration scales]\label{lem:scale-choice}
	Under \eqref{eq:absolute-truncation-condition}, there are $R_n\to\infty$ and $\varepsilon_n\to0$ such that
	\begin{equation}\label{eq:scale-choice-errors}
		R_n^{-\kappa}=o(\varepsilon_n^{\gamma_q}),
		\qquad
		\varepsilon_n^2M(U_{R_n})=o(1).
	\end{equation}
	Independently of that condition, for every
	\begin{equation}\label{eq:relative-eta-range}
		0<\eta<\min\{\kappa/2,2\},
	\end{equation}
	choose $T_n\to\infty$, $\varepsilon_n=e^{-T_n}$, and
	\begin{equation}\label{eq:relative-scale-choice-proof}
		R_n=e^{3T_n/4}.
	\end{equation}
	Then
	\begin{equation}\label{eq:interaction-small}
		R_n^{-\kappa}=o(e^{-\eta T_n}),\qquad
		e^{-2T_n}M(U_{R_n})=o(e^{-\eta T_n}),\qquad
		(1+T_n)e^{-\kappa T_n/2}=o(e^{-\eta T_n}).
	\end{equation}
	The shifted cutoff is supported in a slab whose left endpoint tends to $+\infty$.
\end{lemma}

\begin{proof}
	To establish the first assertion, suppose that \eqref{eq:absolute-truncation-condition} holds and set $R_n=\varepsilon_n^{-\tau}$ for a parameter $\tau > 0$ to be chosen. We distinguish three cases for $\kappa$:
	
	When $\kappa>2$, the mass $M(U_R)$ is uniformly bounded. Choosing any $\tau>\gamma_q/\kappa$ immediately yields both asymptotic limits in \eqref{eq:scale-choice-errors}.
	
	When $\kappa=2$, we select $\tau \in (\gamma_q/2, 1)$. Since $M(U_R)\leq C(1+\log R)$, it follows that
	\[
	\varepsilon_n^2M(U_{R_n})
	\leq C\varepsilon_n^2(1+\tau|\log\varepsilon_n|)=o(1).
	\]
	
	When $0<\kappa<2$, we choose $\tau$ such that
	\[
	\frac{\gamma_q}{\kappa}<\tau<\frac{2}{2-\kappa}.
	\]
	This choice is admissible since condition \eqref{eq:absolute-truncation-condition} for $\kappa<2$ is equivalent to the strict inequality $\gamma_q(2-\kappa)<2\kappa$. By \eqref{eq:truncated-mass}, we have
	\[
	\varepsilon_n^2M(U_{R_n})
	\leq C\varepsilon_n^{2-\tau(2-\kappa)}
	=o(1).
	\]
	while the lower bound on $\tau$ simultaneously guarantees $R_n^{-\kappa}=\varepsilon_n^{\tau\kappa}=o(\varepsilon_n^{\gamma_q})$.
	
	Turning to the second assertion, the choice of $R_n$ and the restriction \eqref{eq:relative-eta-range} on $\eta$ yield $R_n^{-\kappa}=e^{-3\kappa T_n/4}=o(e^{-\eta T_n})$. By Corollary \ref{cor:truncated-extremals}, we have
	\begin{equation}\label{eq:relative-mass-explicit}
		e^{-2T_n}M(U_{R_n})
		\leq C(1+T_n)e^{-\left[2-\frac34(2-\kappa)_+\right]T_n},
		\qquad \text{where } (2-\kappa)_+=\max\{2-\kappa,0\}.
	\end{equation}
	For $\kappa\geq 2$, the exponent in the upper bound is $2$, which strictly exceeds $\eta$. For $0<\kappa<2$, the exponent evaluates to $1/2+3\kappa/4$, which is strictly greater than $\kappa$ and therefore greater than $\eta$. The third estimate in \eqref{eq:interaction-small} is an immediate consequence of the hypothesis $\eta<\kappa/2$. 
	
	Finally, the support of the shifted cutoff on the cylinder is contained in
	\[
	[T_n/4-\log2,7T_n/4+\log2]\times\mathbb S^{N-1},
	\]
	whose left endpoint tends to $+\infty$ as $n \to \infty$.
\end{proof}
Although the above construction yields test functions with the desired critical behavior, the associated mass vanishes upon translation. To restore the prescribed mass, we select a function $h\in C_c^{\infty}(\mathbb R^N\setminus\{0\})$ normalized such that $M(h)=1$. For $\sigma<0$, we set $h_{\sigma}=\sigma\star h$, which yields the scaling relations
\[
M(h_{\sigma})=1,  \quad A(h_{\sigma})=e^{2\sigma}A(h),  \quad B_q(h_{\sigma})=e^{r_q\sigma}B_q(h), \quad\mathcal D(h_{\sigma})=e^{2p\sigma}\mathcal D(h).
\]
For any target mass $0<m\leq\rho^2$, the function $\sqrt m\,h_{\sigma}$ satisfies $M(\sqrt m\,h_\sigma)=m$. As $\sigma\to-\infty$, the three corresponding energy components tend to zero. Moreover, the support of $h_{\sigma}$ is contained in an annulus whose inner radius tends to infinity as $\sigma\to-\infty$. Consequently, by choosing $\sigma$ sufficiently negative, we can ensure that $\operatorname{supp}h_{\sigma}$ is disjoint from any fixed compactly supported profile.

\subsection{The mass-subcritical regime}\label{sec:subcritical}

\begin{lemma}[Branch separation]\label{lem:branch-separation}
	Under the assumptions of Lemma \ref{lem:subcritical-fibers}, after decreasing $\beta_{\mathrm{fib}}$, there are $C,c_0,c_1,c_2>0$ such that
	\begin{equation}\label{eq:lower-branch-small}
		u\in\mathcal M_{\rho,\beta}^{+}
		\quad\Longrightarrow\quad
		A(u)\leq C\beta^{2/(2-r_q)},
	\end{equation}
	while
	\begin{equation}\label{eq:upper-branch-separated}
		u\in\mathcal M_{\rho,\beta}^{-}
		\quad\Longrightarrow\quad
		A(u)\geq c_0,\qquad E_\beta(u)\geq c_1,
		\qquad Q_\beta(u)\leq-c_2,\qquad \mathcal D(u)\geq\tfrac12 A(u).
	\end{equation}
\end{lemma}

\begin{proof}
	Suppose first that $u\in\mathcal M_{\rho,\beta}^{+}$. The strict positivity $Q_{\beta}(u)>0$ implies that
	\[
	\mathcal D(u)<\frac{2-r_q}{2p-2}\beta\delta_qB_q(u).
	\]
	Combining this bound with the identity $P_{\beta}(u)=0$ and Lemma \ref{lem:weighted-interpolation}, we deduce that
	\[
	A(u)\leq C\beta B_q(u)
	\leq C\beta A(u)^{r_q/2}.
	\]
	Since $r_q<2$, the desired estimate \eqref{eq:lower-branch-small} follows immediately.
	
	Next, we consider $u\in\mathcal M_{\rho,\beta}^{-}$. Since $Q_\beta(u)<0$, we obtain the reverse inequality
	\[
	\mathcal D(u)>\frac{2-r_q}{2p-2}\beta\delta_qB_q(u).
	\]
	Since $P_\beta(u)=0$ ensures $A(u)=\beta\delta_qB_q(u)+\mathcal D(u)$, the above lower bound implies $A(u)\leq C\mathcal D(u)$. Since $\mathcal D(u)\leq S_{\mu}(a,b)^{-p}A(u)^p$, we obtain
	\[
	A(u)\leq C A(u)^p,
	\]
	which forces $A(u)\geq c_0$ for some constant $c_0>0$. On the Pohozaev manifold, the energy is given by
	\begin{equation}\label{eq:energy-on-M-A}
		E_{\beta}(u)=\frac{p-1}{2p}A(u)
		-\beta\left(\frac1q-\frac{\delta_q}{2p}\right)B_q(u).
	\end{equation}
	Using \eqref{eq:Bq-shell} and the assumption $r_q/2<1$, we can bound the energy from below for $A(u)\geq c_0$,
	\[
	E_{\beta}(u)\geq cA(u)-C\beta A(u)^{r_q/2}
	\geq \left(c-C\beta c_0^{r_q/2-1}\right)A(u).
	\]
	By choosing $\beta_{\mathrm{fib}}$ sufficiently small, we can guarantee that this coefficient is strictly positive. Furthermore, the upper bound on $B_q(u)$ yields $\beta\delta_qB_q(u)\leq\beta\delta_qK_qc_0^{r_q/2-1}A(u)$. We can restrict $\beta_{\mathrm{fib}}$ further to satisfy
	\[
	\beta\delta_qK_qc_0^{r_q/2-1}\leq\min\left\{\frac12,\frac{p-1}{2p-r_q}\right\}.
	\]
	This choice enforces $\mathcal D(u)\geq A(u)/2$, thereby implying
	\[
	Q_\beta(u)=-(2p-2)A(u)+(2p-r_q)\beta\delta_qB_q(u)
	\leq-(p-2)A(u)\leq-(p-2)c_0.
	\]
	This strictly negative upper bound confirms \eqref{eq:upper-branch-separated}, thus completing the proof.
\end{proof}

Suppose that $2<q<q_c$ and $0<\beta<\beta_{\mathrm{fib}}$. We choose a constant $R_0 \in (0, \sqrt{c_0})$, where $c_0>0$ is given by Lemma \ref{lem:branch-separation}, sufficiently small such that
\[
\frac12R_0^2-\frac{S_{\mu}(a,b)^{-p}}{2p}R_0^{2p}  \geq\frac38R_0^2.
\]
By decreasing $\beta>0$ further if necessary, we may
also assume that
\[
C\beta^{2/(2-r_q)}<R_0^2,
\qquad
\frac{\beta K_q}{q}R_0^{r_q}
\leq \frac18 R_0^2.
\]
With these choices in place, we introduce the set
\[
G_{\rho}=\left\{u\in S_{\rho} : A(u)^{1/2}<R_0\right\}
\]
and define the corresponding infimum level $m_{\beta}^{+}=\inf_{G_{\rho}}E_{\beta}$.

\begin{lemma}[Local minimization]\label{lem:local-well}
	For small $\beta$,
	\begin{equation}\label{eq:local-well-levels}
		-\infty<m_{\beta}^{+}<0,
		\qquad
		\inf_{\substack{u\in S_\rho\\ A(u)^{1/2}=R_0}}
		E_\beta(u)\geq\frac14R_0^2.
	\end{equation}
	Moreover,
	\begin{equation}\label{eq:branches-in-well}
		\mathcal M_{\rho,\beta}^{+}\subset G_{\rho},
		\qquad
		\mathcal M_{\rho,\beta}^{-}\cap G_{\rho}=\varnothing.
	\end{equation}
\end{lemma}

\begin{proof}
	The explicit lower bound
	\[
	E_{\beta}(u)  \geq\frac12A(u)-\frac{\beta K_q}{q}A(u)^{r_q/2} - \frac{S_{\mu}(a,b)^{-p}}{2p}A(u)^p
	\]
	guarantees that $E_\beta$ is bounded from below on the local domain $G_\rho$ while remaining strictly positive on the boundary $\partial G_\rho$. On the other hand, for a fixed $u\in S_{\rho}$, rescaling by $t$ reveals that
	\[
	E_{\beta}(t\star u)  =-\frac{\beta}{q}e^{r_qt}B_q(u)+o(e^{r_qt})<0
	\]
	whenever $t$ is chosen sufficiently negative. Consequently, the infimum must satisfy $m_{\beta}^{+}<0$. Finally, the branch inclusions in \eqref{eq:branches-in-well} follow immediately from Lemma \ref{lem:branch-separation}.
\end{proof}

\begin{lemma}[A minimizing Palais--Smale sequence]\label{lem:well-PSP}
	There exists $u_n\in G_{\rho}$ such that
	\begin{equation}\label{eq:well-PSP}
		E_{\beta}(u_n)\to m_{\beta}^{+},
		\qquad
		d(E_{\beta}|_{S_{\rho}})(u_n)\to0,
		\qquad P_{\beta}(u_n)\to0.
	\end{equation}
\end{lemma}

\begin{proof}
	Choose $\varepsilon_n\downarrow0$ and $v_n\in G_\rho$
	such that $E_\beta(v_n)\leq m_\beta^++\varepsilon_n^2$.
	We equip $\mathbb R\times S_\rho$ with the product metric
	\[
	d((t,u),(s,v))
	=\bigl(|t-s|^2+\|u-v\|_X^2\bigr)^{1/2}
	\]
	and define
	\[
	\widetilde G_\rho
	=\{(t,u)\in\mathbb R\times S_\rho:
	e^{2t}A(u)<R_0^2\},
	\qquad
	\widetilde E_\beta(t,u)=E_\beta(t\star u).
	\]
	Since $M$ is continuous, $S_\rho$ is closed in $X$.
	Hence $\overline{\widetilde G_\rho}$ is complete under
	the product metric $d$ and
	$\inf_{\overline{\widetilde G_\rho}}\widetilde E_\beta
	=m_\beta^+$. By Ekeland's variational principle, we find
	$y_n=(t_n,\widetilde v_n)\in
	\overline{\widetilde G_\rho}$ such that
	\[
	\widetilde E_\beta(y_n)\leq E_\beta(v_n),\qquad
	d(y_n,(0,v_n))\leq\varepsilon_n,
	\]
	and
	\[
	\widetilde E_\beta(z)\geq
	\widetilde E_\beta(y_n)-\varepsilon_n d(z,y_n)
	\quad\text{for all }z\in\overline{\widetilde G_\rho}.
	\]
	Since $m_\beta^+<0$ and the boundary estimate
	\eqref{eq:local-well-levels} gives a positive lower bound
	for $\widetilde E_\beta$ on $\partial\widetilde G_\rho$,
	we have $y_n\in\widetilde G_\rho$ for all sufficiently large $n$.
	We can then apply the variational inequality along
	curves through $y_n$ in both directions to obtain
	\begin{equation}\label{eq:ekeland-well-closeness-proof}
		|t_n|+\|\widetilde v_n-v_n\|_X\to0,\qquad
		\|d\widetilde E_\beta(y_n)\|\leq\varepsilon_n.
	\end{equation}
	
	Set $u_n=t_n\star\widetilde v_n$. Then $u_n\in G_\rho$
	and $E_\beta(u_n)\to m_\beta^+$. Moreover,
	\[
	|P_\beta(u_n)|
	=|\partial_t\widetilde E_\beta(y_n)|
	\leq\varepsilon_n\to0.
	\]
	Fix $\zeta\in T_{u_n}S_\rho$ and define
	$\xi_n=(-t_n)\star\zeta$. Then
	$\xi_n\in T_{\widetilde v_n}S_\rho$.
	By scaling estimates (\ref{eq:dilation-operator-norms}), we have
	\[
	\|\xi_n\|_X
	\leq\max\{1,e^{-t_n}\}\|\zeta\|_X
	\leq C\|\zeta\|_X.
	\]
Using $t_n\star\xi_n=\zeta$ and the chain rule, we obtain
	\[
		|d(E_\beta|_{S_\rho})(u_n)[\zeta]|
		=|d_u\widetilde E_\beta(y_n)[\xi_n]|
		\leq C\varepsilon_n\|\zeta\|_X.
	\]
	Thus $d(E_\beta|_{S_\rho})(u_n)\to0$. Finally,
	$A(u_n)<R_0^2$ and $M(u_n)=\rho^2$, so $(u_n)$
	is bounded in $X$. This completes the proof.
\end{proof}

\begin{proposition}[Negative-energy normalized ground state]
	\label{prop:negative-ground-state}
	For $\beta>0$ sufficiently small, the level $m_\beta^+$ is attained
	by a positive normalized solution $u_\beta^+$.
	Moreover, $u_\beta^+$ is a local minimizer on $S_\rho$
	and a normalized ground state.
\end{proposition}

\begin{proof}
	Let $\{u_n\} \subset G_\rho$ be the sequence provided by Lemma \ref{lem:well-PSP}. This sequence is bounded in $X$ because $A(u_n) < R_0^2$ and $M(u_n) = \rho^2$, which allows us to apply Theorem \ref{thm:global-compactness-intro}\textit{(ii)}. Since the strict inequality $m_\beta^+ < 0$ precludes a trivial weak limit, we can pass to a subsequence,
	\[
	u_n\rightharpoonup u^0 \quad\text{in }X,
	\qquad
	u^0 \not\equiv 0,
	\qquad
	M(u^0)=\rho^2.
	\]
Combining weak lower semicontinuity with the energy
decomposition, we obtain
	\[
	A(u^0)\leq R_0^2,
	\qquad
	E_\beta(u^0)
	=m_\beta^+-\sum_{j=1}^{k}J(U^j)
	\leq m_\beta^+<0.
	\]
	If we assume $A(u^0)=R_0^2$, then \eqref{eq:local-well-levels}
	would imply
	\[
	E_\beta(u^0)\geq\frac14R_0^2>0,
	\]
	which is a contradiction.
Thus $u^0\in G_\rho$, so $E_\beta(u^0)\geq m_\beta^+$.
	It follows that
	\[
	m_\beta^+
	=E_\beta(u^0)+\sum_{j=1}^{k}J(U^j)
	\geq m_\beta^++kT_{\mu,a,b}.
	\]
	Thus $k=0$. Applying Theorem
	\ref{thm:global-compactness-intro}\textit{(ii)}, we obtain
	\[
	u_n\to u^0 \quad\text{in }X,
	\qquad
	E_\beta(u^0)=m_\beta^+.
	\]
	
	Define $u_\beta^+=|u^0|$.
	Since the functionals $M$, $A$, and $E_\beta$ remain invariant under
	the mapping $u\mapsto|u|$, we deduce
	\[
	u_\beta^+\in G_\rho,
	\qquad
	E_\beta(u_\beta^+)=m_\beta^+
	=\inf_{v\in G_\rho}E_\beta(v).
	\]
Since $G_\rho$ is relatively open in $S_\rho$,
we conclude that $u_\beta^+$ is a local minimizer
of $E_\beta$ on $S_\rho$.
	In particular,
	\[
	d(E_\beta|_{S_\rho})(u_\beta^+)=0.
	\]
	Local regularity and the strong maximum principle imply
	$u_\beta^+>0$ in $\mathbb R^N\setminus\{0\}$.
	Moreover, we can apply Proposition \ref{prop:pohozaev} to obtain
	\[
	P_\beta(u_\beta^+)=0,
	\qquad
	\lambda_\beta^+\rho^2
	=-\beta(1-\delta_q)B_q(u_\beta^+)<0.
	\]
	Since $u_\beta^+$ is a local minimizer on $S_\rho$,
	the fiber map $t\mapsto E_\beta(t\star u_\beta^+)$
	has a local minimum at $t=0$. Hence
	$Q_\beta(u_\beta^+)\geq0$. Since
	$\mathcal M_{\rho,\beta}^{0}=\varnothing$, we conclude that
	$u_\beta^+\in\mathcal M_{\rho,\beta}^{+}$.
	
	Finally, let $v\in S_\rho$ be an arbitrary normalized solution.
	By Proposition \ref{prop:pohozaev} and
	Lemma \ref{lem:subcritical-fibers}, we have
	\[
	v\in\mathcal M_{\rho,\beta}
	=\mathcal M_{\rho,\beta}^+
	\cup\mathcal M_{\rho,\beta}^-.
	\]
	If $v\in\mathcal M_{\rho,\beta}^+$, we use
	\eqref{eq:branches-in-well} to conclude that $v\in G_\rho$, thus
	$E_\beta(v)\geq m_\beta^+$.
	If $v\in\mathcal M_{\rho,\beta}^-$, we use
	Lemma \ref{lem:branch-separation} to obtain
	$E_\beta(v)>0>m_\beta^+$.
	Therefore, $u_\beta^+$ minimizes the energy among all
	normalized solutions and is indeed a normalized ground state.
\end{proof}

We define the infimum of the energy on the upper branch as
\[
m_{\beta}^{-}=\inf_{u\in\mathcal M_{\rho,\beta}^{-}}E_{\beta}(u).
\]
By the fiber geometry established in Lemma \ref{lem:subcritical-fibers}, this infimum admits the equivalent minimax characterization
\[
m_{\beta}^{-}  =\inf_{v\in S_{\rho}}\max_{t\in\mathbb R}E_{\beta}(t\star v).
\]
Consequently, $m_{\beta}^{-}$ corresponds to a fiber mountain-pass level.

\begin{lemma}[Upper-branch Palais--Smale sequence]
	\label{lem:upper-PS}
	There exists a bounded sequence
	$\{u_n\}\subset\mathcal M_{\rho,\beta}^-$ such that
	\[
	E_\beta(u_n)\to m_\beta^-,
	\qquad
	\|d(E_\beta|_{S_\rho})(u_n)\|\to0.
	\]
\end{lemma}

\begin{proof}
Let $\widehat u_n\in\mathcal M_{\rho,\beta}^-$ be a minimizing sequence.
By \eqref{eq:energy-on-M-A} and \eqref{eq:Bq-shell}, we have
\[
E_\beta(\widehat u_n)\geq\frac{p-1}{2p}A(\widehat u_n)
-C\beta A(\widehat u_n)^{r_q/2}.
\]
Since $(\widehat u_n)$ is a minimizing sequence,
the energies $E_\beta(\widehat u_n)$ are bounded.
Because $r_q/2<1$, the preceding estimate implies
that $A(\widehat u_n)$ is bounded above. Combining
this bound with Lemma \ref{lem:branch-separation},
we obtain
\[
0<c\leq A(\widehat u_n)\leq C,
\qquad
Q_\beta(\widehat u_n)\leq-c.
\]
Set $s_n=-\tfrac12\log A(\widehat u_n)$ and
$\widehat v_n=s_n\star\widehat u_n$. Then
$A(\widehat v_n)=1$, $|s_n|\leq C$ and
$t_2(\widehat v_n)=-s_n$. Thus $\widehat v_n$ is a minimizing
sequence for $\Psi_2$ on $S_\rho$.
Since $\Psi_2\in C^1(S_\rho)$ and $S_\rho$ is complete,
we apply Ekeland's variational principle to find
$v_n\in S_\rho$ such that
\[
\|v_n-\widehat v_n\|_X\to0,\qquad
\Psi_2(v_n)\to m_\beta^-,\qquad \|d\Psi_2(v_n)\|\to0.
\]
Set $t_n=t_2(v_n)$ and $u_n=t_n\star v_n$.
By Lemma \ref{lem:uniform-projection}, we have
$|t_n|\leq C$.
For any $\zeta\in T_{u_n}S_\rho$, combining
\eqref{eq:reduced-derivative-formula} with
\eqref{eq:dilation-operator-norms}, we obtain
\[
\bigl|d(E_\beta|_{S_\rho})(u_n)[\zeta]\bigr|
=\bigl|d\Psi_2(v_n)[(-t_n)\star\zeta]\bigr|
\leq C\|d\Psi_2(v_n)\|\,\|\zeta\|_X.
\]
Hence the constrained derivative tends to zero, while
$E_\beta(u_n)=\Psi_2(v_n)\to m_\beta^-$.
Moreover,
since $\|v_n-\widehat v_n\|_X\to0$ and
$A(\widehat v_n)=1$, we have $A(v_n)\to1$.
Since $(t_n)$ is bounded, we obtain
\[\|u_n\|_X^2=e^{2t_n}A(v_n)+\rho^2\leq C.\]
\end{proof}

\begin{lemma}[Decay of the local minimizer]\label{lem:lower-tail}
	There exists $C>0$ such that
	\[
	w_{\beta}^{+}(t,\omega)
	+|\partial_tw_{\beta}^{+}(t,\omega)|
	+|\nabla_{\omega}w_{\beta}^{+}(t,\omega)|
	\leq Ce^{-\kappa t/2}
	\quad\text{for }t\geq1.
	\]
\end{lemma}

\begin{proof}
Set $w=w_\beta^+$, $\nu=\kappa/2$, and $V=\mathcal K_\mu*w^p$.
Then
\[
L_\kappa w=\lambda_\beta^+e^{-2t}w+F_+,
\qquad F_+=\beta e^{-\gamma_qt}w^{q-1}+Vw^{p-1}.
\]
By Lemma \ref{lem:full-equation-regularity}, $V$ is bounded on the
right half-cylinder and $w(t,\cdot)\to0$ uniformly.
Since $\lambda_\beta^+<0$ and $p,q>2$, we have
$L_\kappa w\leq F_+\leq\varepsilon w$ for large $t$.
The comparison argument in Proposition \ref{prop:axial-decay} gives
$w(t,\omega)\leq C_\sigma e^{-\sigma t}$ for every $0<\sigma<\nu$.
Choose $\sigma$ close enough to $\nu$ that
\[
b_2=\min\{\gamma_q+(q-1)\sigma,(p-1)\sigma\}>\nu.
\]
Then $F_+\leq Ce^{-b_2t}$. Fix $T\geq1$ and $C_0>0$ such that
$F_+\leq C_0e^{-b_2(t-T)}$ for $t\geq T$, and define
\[
\Phi(t)=Ae^{-\nu(t-T)}
+\frac{C_0}{b_2^2-\nu^2}
\bigl(e^{-\nu(t-T)}-e^{-b_2(t-T)}\bigr),
\qquad A\geq\|w(T,\cdot)\|_\infty.
\]
Since $w,\Phi\in H^1(\mathcal C_T)$ and
$w(T,\cdot)\leq\Phi(T)$, we have
$(w-\Phi)_+\in H_0^1(\mathcal C_T)$.
We can therefore apply the comparison argument
from Proposition \ref{prop:axial-decay} to conclude that
$w\leq\Phi\leq Ce^{-\nu t}$ for $t\geq T$.

Using the decay estimates for $w$ and $F=L_\kappa w$, we obtain
\[
\|w\|_{L^\infty(K_2(t))}
+\|L_\kappa w\|_{L^\infty(K_2(t))}
\leq Ce^{-\nu t},
\qquad t\geq T+2.
\]
Applying the interior elliptic estimate from Proposition
\ref{prop:axial-decay}, we obtain the same decay bound
for $|\partial_t w|+|\nabla_\omega w|$.
We use local regularity to obtain the same estimate on the
remaining compact slab $1\leq t\leq T+2$.
\end{proof}

\begin{lemma}[Estimates for the test path]\label{lem:insertion-path}
	Let $u_\beta^+$ be the local minimizer given by Proposition \ref{prop:negative-ground-state}, with $\beta>0$ sufficiently small, and write $w=w_\beta^+$. Fix $\eta$ satisfying \eqref{eq:relative-eta-range}, choose the sequences in \eqref{eq:relative-scale-choice-proof}, and put
	\[
	Z_n(t,\omega)=W_{R_n}(t-T_n,\omega),\qquad
	m_n=\int_{\mathcal C}e^{-2t}Z_n^2\,dt\,d\omega,\qquad
	c_n=\int_{\mathcal C}e^{-2t}wZ_n\,dt\,d\omega,
	\]
	where $W_R$ is the cylinder transform of the function $U_R$ normalized by \eqref{eq:UR-A-D}. For every fixed $L>1$ and all sufficiently large $n$,
	there is a positive $C^1$ function
	$\alpha_n:[0,L]\to\mathbb R$ such that
	\begin{equation}\label{eq:mass-corrected-path-cylinder}
		\Gamma_n(s)=\alpha_n(s)w+sZ_n\in S_\rho,
		\quad \|\alpha_n-1\|_{C^1([0,L])}=o(e^{-\eta T_n}).
	\end{equation}
	Define
	\[
	C_n=\mathfrak B_{\mathcal C}(w^p,Z_n^p),
	\qquad g(s)=\frac{s^2}{2}-\frac{s^{2p}}{2p}.
	\]
	There exists $\mathcal R_n\geq0$, independent of $s\in[0,L]$, such that
	\begin{equation}\label{eq:relative-error-explicit}
			\frac{\mathcal R_n}{e^{-\eta T_n}}\leq C_{\beta,L}\Bigl[
			e^{-(3\kappa/4-\eta)T_n}+(1+T_n)e^{-(\kappa/2-\eta)T_n}
			+(1+T_n)e^{-[2-\frac34(2-\kappa)_+-\eta]T_n}\Bigr]\longrightarrow0.
	\end{equation}
	Then $C_n\geq0$, and uniformly on $[0,L]$,
	\begin{equation}\label{eq:relative-path-energy}
			E_\beta(\Gamma_n(s))\leq{}m_\beta^++A(U_{R_n})g(s)
			-\frac\beta q s^q e^{-\gamma_qT_n}B_q(U_{R_n})
			-\frac{\alpha_n(s)^p s^p}{p}C_n+\mathcal R_n,
	\end{equation}
	and
	\begin{equation}\label{eq:uniform-path-P}
		\sup_{0\leq s\leq L}
		\left|P_\beta(\Gamma_n(s))-A(U)(s^2-s^{2p})\right|\to0.
	\end{equation}
	For $L=2$, the path intersects $\mathcal M_{\rho,\beta}^-$.
\end{lemma}

\begin{proof}
	Recalling the definition $Z_n(t,\omega)=W_{R_n}(t-T_n,\omega)$, we deduce from the decay estimate for $W$ and the cutoff construction that
	\[
	|Z_n|+|\partial_tZ_n|+|\nabla_\omega Z_n|
	\leq
	Ce^{-\kappa|t-T_n|/2}
	\]
	on the support of $Z_n$. Moreover, for $n$ sufficiently large, we have $\operatorname{supp}Z_n\subset\{t>1\}$.
	
	We use Lemma \ref{lem:lower-tail} to estimate the interaction
	terms and obtain
	\[
	\begin{aligned}
		|\langle w,Z_n\rangle_\kappa|
		\leq
		C\int_0^\infty
		e^{-\kappa t/2}
		e^{-\kappa|t-T_n|/2}\,dt
		&=
		C\int_0^{T_n}
		e^{-\kappa T_n/2}\,dt
		+
		C\int_{T_n}^{\infty}
		e^{-\kappa t/2}
		e^{-\kappa(t-T_n)/2}\,dt\\
		&\leq
		C\left(T_n+\frac1\kappa\right)
		e^{-\kappa T_n/2}.
	\end{aligned}
	\]
	Since the derivatives of the cutoff function are uniformly bounded, the terms involving these derivatives satisfy the same estimate. Consequently,
	\begin{equation}\label{eq:interaction-estimates}
		|\langle w,Z_n\rangle_\kappa|+c_n
		\leq
		C\left(T_n+\frac1\kappa\right)
		e^{-\kappa T_n/2}.
	\end{equation}
	Combining this bound with \eqref{eq:interaction-small}, we have
	\[
	m_n+c_n=o(e^{-\eta T_n}).
	\]
	
To preserve the mass constraint, $\alpha_n(s)$ must satisfy
\[
\rho^2\alpha_n(s)^2
+2sc_n\alpha_n(s)+s^2m_n=\rho^2.
\]
For each $s\in[0,L]$, define
\[
f_n(\alpha)
=\rho^2\alpha^2+2sc_n\alpha+s^2m_n-\rho^2.
\]
Since $m_n+c_n\to0$, we have
$f_n(1/2)<0\leq f_n(1)$ for all $s\in[0,L]$ and
all sufficiently large $n$.
Moreover,
\[
f_n'(\alpha)=2\rho^2\alpha+2sc_n>0
\qquad\text{for }\alpha>0.
\]
Thus there is a unique positive root
$\alpha_n(s)\in(1/2,1]$, with $\alpha_n(0)=1$.
The implicit function theorem gives
$\alpha_n\in C^1([0,L])$.
We rearrange the mass constraint to obtain
\[
1-\alpha_n(s)
=
\frac{2\alpha_n(s)sc_n+s^2m_n}
{\rho^2(1+\alpha_n(s))}
\leq C_L(c_n+m_n).
\]
Differentiating the mass constraint with respect to $s$, we obtain
\[
\alpha_n'(s)
=
-\frac{\alpha_n(s)c_n+sm_n}
{\rho^2\alpha_n(s)+sc_n}.
\]
	Using $w\geq0$ and $Z_n\geq0$, we obtain
	\[
	(\alpha_n(s)w+sZ_n)^p
	\geq
	\alpha_n(s)^pw^p+s^pZ_n^p .
	\]
We use the positivity of the kernel to obtain
	\begin{equation}\label{eq:retained-nonlocal-interaction}
		\begin{aligned}
			\mathcal D_{\mathcal C}(\Gamma_n(s))
			&\geq\iint_{\mathcal C\times\mathcal C}
			\mathcal K_\mu(t-r,\omega,\omega')
			\bigl[\alpha_n(s)^pw(t,\omega)^p+s^pZ_n(t,\omega)^p\bigr]\\
			&\qquad\times\bigl[\alpha_n(s)^pw(r,\omega')^p+s^pZ_n(r,\omega')^p\bigr]
			\,dt\,d\omega\,dr\,d\omega'\\
			&=\alpha_n(s)^{2p}\mathcal D_{\mathcal C}(w)
			+s^{2p}\mathcal D_{\mathcal C}(Z_n)
			+2\alpha_n(s)^ps^pC_n,
		\end{aligned}
	\end{equation}
	where
	\[
	C_n=\iint_{\mathcal C\times\mathcal C}
	\mathcal K_\mu(t-r,\omega,\omega')
	w(t,\omega)^pZ_n(r,\omega')^p\,dt\,d\omega\,dr\,d\omega'.
	\]
	A similar argument yields the lower bound for the local nonlinear term
	\[
	B_q(\Gamma_n(s))
	\geq
	\alpha_n(s)^qB_q(w)+s^qB_q(Z_n).
	\]
	Meanwhile, the quadratic part expands as
	\[
	A(\Gamma_n(s))
	=
	A(\alpha_n(s)w)
	+s^2A(Z_n)
	+2\alpha_n(s)s\langle w,Z_n\rangle_\kappa .
	\]
	Combining the Lipschitz estimate $|E_\beta(\alpha_n(s)w)-E_\beta(w)|\leq C_L|\alpha_n(s)-1|$, the asymptotic expansion $A(U_{R_n})=A(U)+O(R_n^{-\kappa})$ and the interaction bounds in \eqref{eq:interaction-estimates}, we see that all resulting error terms are bounded by
	\[
	\mathcal R_n := C_{\beta,L}
	\left(
	R_n^{-\kappa}
	+(1+T_n)e^{-\kappa T_n/2}
	+e^{-2T_n}M(U_{R_n})
	\right).
	\]
	Together with the identities $A(U_{R_n})=\mathcal D(U_{R_n})$ and $B_q(Z_n)=e^{-\gamma_qT_n}B_q(U_{R_n})$, this establishes \eqref{eq:relative-path-energy}.
	
	It remains to verify \eqref{eq:uniform-path-P}. For any test function $\varphi\in C_c^\infty(\mathcal C)$, we have
	\[
	|\langle Z_n,\varphi\rangle_\kappa|
	\leq
	\|W_{R_n}-W\|_\kappa\|\varphi\|_\kappa
	+
	|\langle W(\cdot-T_n,\cdot),\varphi\rangle_\kappa|.
	\]
	As $n\to\infty$, the first term on the right-hand side vanishes because $W_{R_n}\to W$ strongly in $H^1(\mathcal C)$, while the second term tends to zero since $T_n\to\infty$ and $\varphi$ is compactly supported. By density, it follows that
	\[
	Z_n\rightharpoonup0
	\quad\text{in }H^1(\mathcal C).
	\]
	Since $m_n\to0$, we apply the weighted interpolation inequality to obtain
	\[
	B_q(Z_n)
	\leq
	C A(Z_n)^{r_q/2}m_n^{\gamma_q/2}
	\to0 .
	\]
	Moreover,
	\[
	\sup_{s\in[0,L]}B_q(\Gamma_n(s)-w)^{1/q}
	\leq
	\|\alpha_n-1\|_\infty B_q(w)^{1/q}
	+L B_q(Z_n)^{1/q}\to0.
	\]
	Hence $B_q(\Gamma_n(s))=B_q(w)+o(1)$ uniformly on $[0,L]$. 
	
	Since the supports of $Z_n$ move to $+\infty$, we have
	$Z_n\to0$ pointwise a.e. on $\mathcal C$. We can therefore
	apply the nonlocal Brézis--Lieb splitting to obtain, for each
	fixed $s\geq0$,
	\[
	\mathcal D_{\mathcal C}(w+sZ_n)
	=
	\mathcal D_{\mathcal C}(w)
	+s^{2p}\mathcal D_{\mathcal C}(Z_n)
	+o(1).
	\]
	Since the derivative bound $\|\mathcal D_{\mathcal C}'(v)\|_{H^{-1}}\leq C\|v\|_{H^1}^{2p-1}$ ensures that the family of functions $s\mapsto \mathcal D_{\mathcal C}(w+sZ_n)$ is equi-Lipschitz on bounded intervals, a standard finite covering argument upgrades this pointwise limit to uniform convergence on $[0,L]$. Moreover, by the estimate on $|\alpha_n(s)-1|$, replacing $w$ with $\alpha_n(s)w$ introduces only an error of order $o(1)$ uniformly in $s\in[0,L]$. Combining these facts with \eqref{eq:interaction-estimates}, the identity $P_\beta(w)=0$ and the convergence $A(U_{R_n})\to A(U)$, we obtain \eqref{eq:uniform-path-P}.
	
	Finally, we take $L=2$ and $s_0=1/2$. The limiting function
	$A(U)(s^2-s^{2p})$ in \eqref{eq:uniform-path-P} is positive
	at $s=s_0$ and negative at $s=2$. By uniform convergence,
	for all sufficiently large $n$ we can apply the intermediate
	value theorem to find $s_n^*\in(s_0,2)$ such that
	\[
	P_\beta(\Gamma_n(s_n^*))=0.
	\]
	Using \eqref{eq:interaction-estimates}, we obtain
	\[
		A(\Gamma_n(s_n^*))
		=A(\alpha_n(s_n^*)w)
		+(s_n^*)^2A(U_{R_n})+o(1)
		\geq s_0^2A(U_{R_n})+o(1)
		\geq \frac12s_0^2A(U)>0
	\]
	for all sufficiently large $n$. Fixing $\beta>0$ sufficiently small so that the lower-branch threshold $C\beta^{2/(2-r_q)}$ from Lemma \ref{lem:branch-separation} is strictly smaller than $\frac12s_0^2A(U)$, we deduce that $\Gamma_n(s_n^*)\notin\mathcal M_{\rho,\beta}^+$. Since $\mathcal M_{\rho,\beta}^0=\varnothing$, we conclude that
	\[
	\Gamma_n(s_n^*)\in\mathcal M_{\rho,\beta}^{-}
	\]
	for all sufficiently large $n$.
\end{proof}

\begin{proposition}[Relative strict level]\label{prop:relative-strict-level}
Assume $2<q<q_c$, \eqref{eq:decay-range}, and either
$\gamma_q<\kappa/2$ or $\mu<\kappa$.
For $\beta>0$ sufficiently small,
\begin{equation}\label{eq:relative-strict-level}
m_\beta^-<m_\beta^++T_{\mu,a,b}.
\end{equation}
\end{proposition}

\begin{proof}
	Fix $\beta>0$ sufficiently small. In the first case, set
	$\eta=\gamma_q$; in the second, set $\eta=\mu/2$. In either
	case, $0<\eta<\min\{\kappa/2,2\}$, so we can apply
	Lemma \ref{lem:insertion-path} with $L=2$ and
	$R_n=e^{3T_n/4}$. In particular, we have
	\[
	\mathcal R_n=o(e^{-\eta T_n}),
	\qquad R_n^{-\kappa}=o(e^{-\eta T_n}).
	\]
	Fix $\delta\in(0,1/2)$ and set $I_\delta=[1-\delta,1+\delta]$.
	We first show that the negative terms in
	\eqref{eq:relative-path-energy} lower the energy by an amount
	of order $e^{-\eta T_n}$, uniformly on $I_\delta$.
	
	If $\gamma_q<\kappa/2$, we use
	\eqref{eq:truncated-B-lower} to obtain
	\[
	\frac{\beta}{q}s^qe^{-\gamma_qT_n}B_q(U_{R_n})
	\geq c_\beta e^{-\gamma_qT_n},
	\qquad s\in I_\delta.
	\]
	
	If $\mu<\kappa$, let $w=w_\beta^+$, let $W$ denote the
	cylinder transform of $U$, and set
	$K=[-1,1]\times\mathbb S^{N-1}$ and
	$K_n=[T_n-1,T_n+1]\times\mathbb S^{N-1}$.
	For $|t|,|r|\leq1$ and all angular variables, we have
	\[
	2\bigl(\cosh(t-T_n-r)-\omega\cdot\omega'\bigr)
	\leq Ce^{T_n},
	\qquad
	\mathcal K_\mu(t-T_n-r,\omega,\omega')
	\geq ce^{-\mu T_n/2}.
	\]
	For large $n$, we have $W_{R_n}=c_{R_n}W$ on $K$,
	where $c_{R_n}\to1$. Restricting $C_n$ to $K\times K_n$,
	we obtain
	\[
		C_n
		\geq ce^{-\mu T_n/2}
		\left(\int_K w(t,\omega)^p\,dt\,d\omega\right)
		\left(\int_K W_{R_n}(r,\omega')^p\,dr\,d\omega'\right)
		\geq c_\beta e^{-\mu T_n/2}.
	\]
	The last inequality follows from the positivity of $w$ and $W$.
	Since $\alpha_n\to1$ uniformly on $[0,2]$, we further obtain
	\[
	\frac{\alpha_n(s)^ps^p}{p}C_n
	\geq c_\beta'e^{-\mu T_n/2},
	\qquad s\in I_\delta.
	\]
	
	Using $g(s)\leq g(1)$ and
	$A(U_{R_n})g(1)=T_{\mu,a,b}+O(R_n^{-\kappa})$,
	we deduce from \eqref{eq:relative-path-energy} that
	\[
		\max_{s\in I_\delta}E_\beta(\Gamma_n(s))
		\leq m_\beta^++T_{\mu,a,b}
		-c_\beta''e^{-\eta T_n}+o(e^{-\eta T_n})
		<m_\beta^++T_{\mu,a,b}
	\]
	for all sufficiently large $n$.
	
	It remains to estimate the energy outside $I_\delta$.
	Since $g$ has a strict maximum at $s=1$, we can choose
	$d_\delta>0$ such that
	$g(s)\leq g(1)-d_\delta$ for
	$s\in[0,2]\setminus I_\delta$. Dropping the two nonpositive
	terms in \eqref{eq:relative-path-energy}, we obtain
	\[
		\sup_{s\in[0,2]\setminus I_\delta}E_\beta(\Gamma_n(s))
		\leq m_\beta^++T_{\mu,a,b}
		-d_\delta A(U)+o(1)
		<m_\beta^++T_{\mu,a,b}
	\]
	for all sufficiently large $n$.
	
By Lemma \ref{lem:insertion-path}, we can choose
$s_n^*\in[0,2]$ such that
$\Gamma_n(s_n^*)\in\mathcal M_{\rho,\beta}^-$.
The preceding estimates then yield
\eqref{eq:relative-strict-level}.
\end{proof}

We now use the strict energy inequality to prove convergence of the Palais--Smale sequence at $m_\beta^-$.

\begin{proposition}[Second normalized solution]
	\label{prop:second-solution}
	Assume \eqref{eq:main-multiplicity-range}. For $\beta>0$ sufficiently small,
	the level $m_\beta^-$ is attained by a positive
	normalized solution $u_\beta^-$, and
	\[
	0<m_\beta^-<m_\beta^++T_{\mu,a,b}.
	\]
\end{proposition}

\begin{proof}
	By Lemma \ref{lem:branch-separation}, we have
	$m_\beta^->0$. Under \eqref{eq:main-multiplicity-range},
	we apply Proposition \ref{prop:relative-strict-level} to obtain
	\begin{equation}\label{eq:relative-level-used-in-compactness-proof}
		0<m_\beta^-<m_\beta^++T_{\mu,a,b}<T_{\mu,a,b},
	\end{equation}
	where the last inequality follows from $m_\beta^+<0$.
	
	Let $\{u_n\}\subset\mathcal M_{\rho,\beta}^-$ be the bounded
	Palais--Smale sequence provided by Lemma \ref{lem:upper-PS}.
	Since $P_\beta(u_n)=0$, we can apply
	Theorem \ref{thm:global-compactness-intro}\textit{(ii)}. Let $u^0$ denote its weak limit.
	If $u^0 \equiv 0$, the same theorem implies either $m_\beta^-=0$ or $m_\beta^-\geq T_{\mu,a,b}$, contradicting \eqref{eq:relative-level-used-in-compactness-proof}.
	Therefore $u^0 \not\equiv 0$ and Theorem \ref{thm:global-compactness-intro}\textit{(ii)} ensures that $u^0\in S_\rho$ is a normalized solution.
	
	By Proposition \ref{prop:negative-ground-state}, we have
	\[
	E_\beta(u^0)\geq m_\beta^+.
	\]
Using the energy decomposition, we obtain
	\[
	m_\beta^-
	=E_\beta(u^0)+\sum_{j=1}^{k}J(U^j)
	\geq m_\beta^++kT_{\mu,a,b}.
	\]
The strict inequality in
\eqref{eq:relative-level-used-in-compactness-proof}
forces $k=0$. The global compactness theorem then gives
\[
u_n\to u^0\quad\text{in }X,
\qquad
E_\beta(u^0)=m_\beta^-,
\qquad
P_\beta(u^0)=0.
\]
	Moreover, \eqref{eq:upper-branch-separated} and the continuity
	of $Q_\beta$ imply
	\[
	Q_\beta(u^0)=\lim_{n\to\infty}Q_\beta(u_n)\leq-c_2<0.
	\]
	Thus $u^0\in\mathcal M_{\rho,\beta}^-$.
	
	Define $v=|u^0|$. Since the functionals $M$, $E_\beta$, $P_\beta$, and
	$Q_\beta$ remain invariant under the mapping $u\mapsto|u|$, we deduce
	$v\in\mathcal M_{\rho,\beta}^-$ and $t_2(v)=0$. Consequently,
	\[
	\Psi_2(v)=E_\beta(v)=m_\beta^-
	=\inf_{z\in S_\rho}\Psi_2(z).
	\]
	It follows that $d\Psi_2(v)=0$ and
	Lemma \ref{lem:natural-constraint} shows that $v$ is a
	constrained critical point of $E_\beta$ on $S_\rho$.
	We use local regularity and the strong maximum principle to show that
	$v>0$ in $\mathbb R^N\setminus\{0\}$.
	Setting $u_\beta^-=v$ completes the proof.
	The corresponding Lagrange multiplier is strictly negative
	by Proposition \ref{prop:pohozaev}.
\end{proof}

\subsection{The mass-critical and mass-supercritical regimes}\label{sec:critical-supercritical}

Assume $q_c\leq q<2^{\sharp}$.  In the critical case, let $0<\beta<\beta_*$.  Define
\[
m_{\beta}=\inf_{u\in\mathcal M_{\rho,\beta}}E_{\beta}(u).
\]

\begin{lemma}[Coercivity on the Pohozaev manifold]
	\label{lem:coercivity-M}
	The level $m_\beta$ is strictly positive, and every energy-bounded sequence in
	$\mathcal M_{\rho,\beta}$ is bounded in $X$.
	Moreover, there exist constants $c_2,c_3>0$ such that
	\begin{equation}\label{eq:A-Q-uniform}
		A(u)\geq c_2,
		\qquad
		Q_\beta(u)\leq-c_3,
		\qquad \text{for all } u\in\mathcal M_{\rho,\beta}.
	\end{equation}
\end{lemma}

\begin{proof}
	Let $u\in\mathcal M_{\rho,\beta}$, we have
	\[
	A(u)=\beta\delta_qB_q(u)+\mathcal D(u).
	\]
	
	Assume first that $q=q_c$, which implies $r_q=2$ and
	$\delta_q=2/q_c$. By
	\eqref{eq:critical-fiber-positive-coefficient}, there exists
	a constant $c>0$ independent of $u$ such that
	\begin{equation}\label{eq:D-controls-A-critical}
		\mathcal D(u)
		=A(u)-\beta\delta_qB_q(u)\geq cA(u).
	\end{equation}
	Combining this lower bound with
	\[
	\mathcal D(u)\leq S_\mu(a,b)^{-p}A(u)^p,
	\]
	we deduce $A(u)\geq c_2>0$. Moreover, we have
	\[
	\begin{aligned}
		E_\beta(u)
		&=\frac{p-1}{2p}\mathcal D(u)
		\geq \frac{c(p-1)}{2p}A(u),\\
		-Q_\beta(u)
		&=(2p-2)\mathcal D(u)
		\geq c(2p-2)A(u).
	\end{aligned}
	\]
	
	Assume next that $q>q_c$, which yields $r_q>2$.
	By \eqref{eq:Bq-shell} and the definition of $S_\mu(a,b)$, we have
	\[
	A(u)
	\leq \beta\delta_qK_qA(u)^{r_q/2}
	+S_\mu(a,b)^{-p}A(u)^p.
	\]
	Since $A(u)>0$, we have
	\[
	1\leq \beta\delta_qK_qA(u)^{(r_q-2)/2}
	+S_\mu(a,b)^{-p}A(u)^{p-1}.
	\]
	Because both exponents on the right-hand side are strictly positive,
	it follows that $A(u)\geq c_2>0$ uniformly on
	$\mathcal M_{\rho,\beta}$.
	
	Using $r_q=q\delta_q$ and $P_\beta(u)=0$, we deduce
	\begin{equation}\label{eq:energy-super-M}
			E_\beta(u)=\frac{r_q-2}{2r_q}\beta\delta_qB_q(u)
			+\frac{p-1}{2p}\mathcal D(u)
			\geq
			\min\left\{\frac{r_q-2}{2r_q},\frac{p-1}{2p}\right\}A(u).
	\end{equation}
	Similarly,
	\[
-Q_\beta(u)
		=(r_q-2)\beta\delta_qB_q(u)
		+(2p-2)\mathcal D(u)
\geq \min\{r_q-2,2p-2\}A(u).
\]
	
In both cases we can choose constants $c_E,c_Q>0$
independent of $u$ so that
\[
E_\beta(u)\geq c_EA(u)\geq c_Ec_2,
\qquad
-Q_\beta(u)\geq c_QA(u)\geq c_Qc_2.
\]
	Thus $m_\beta\geq c_Ec_2>0$ and \eqref{eq:A-Q-uniform}
	holds with $c_3=c_Qc_2$. Moreover,
	\[
	\|u\|_X^2=A(u)+\rho^2
	\leq c_E^{-1}E_\beta(u)+\rho^2.
	\]
	Hence every sequence in $\mathcal M_{\rho,\beta}$ with bounded
	energy is bounded in $X$.
\end{proof}

\begin{lemma}[Palais--Smale sequence at $m_\beta$]
	\label{lem:m-PS}
	There exists a bounded sequence
	$\{u_n\}\subset\mathcal M_{\rho,\beta}$ such that
	\[
	E_\beta(u_n)\to m_\beta,
	\qquad
	\|d(E_\beta|_{S_\rho})(u_n)\|\to0.
	\]
\end{lemma}

\begin{proof}
Let $\widehat u_n\in\mathcal M_{\rho,\beta}$ be a minimizing
sequence for $m_\beta$. By Lemma \ref{lem:coercivity-M}, we have
\[
0<c\leq A(\widehat u_n)\leq C,
\qquad Q_\beta(\widehat u_n)\leq-c.
\]
The construction in the proof of Lemma \ref{lem:upper-PS} now applies
with the unique projection $t$ and reduced functional $\Psi$ in place
of $t_2$ and $\Psi_2$. Indeed, $\Psi\in C^1(S_\rho)$,
$\inf_{S_\rho}\Psi=m_\beta$ and the preceding bounds verify the
assumptions of Lemma \ref{lem:uniform-projection}. Normalizing
$A(\widehat v_n)=1$ and applying Ekeland's principle therefore gives
$u_n=t(v_n)\star v_n\in\mathcal M_{\rho,\beta}$ with
\[
E_\beta(u_n)\to m_\beta,\qquad
\|d(E_\beta|_{S_\rho})(u_n)\|
\leq C\|d\Psi(v_n)\|\to0.
\]
The boundedness follows from Lemma \ref{lem:coercivity-M}.
\end{proof}

We now prove the required strict energy bound.

\begin{proposition}[Absolute strict level]
	\label{prop:absolute-strict-level}
	Assume $q_c\leq q<2^\sharp$, with $0<\beta<\beta_*$
	when $q=q_c$, and assume
	\eqref{eq:absolute-truncation-condition}. Then
	\begin{equation}\label{eq:absolute-strict-level}
		m_\beta<T_{\mu,a,b}.
	\end{equation}
\end{proposition}

\begin{proof}
	Let $R_n$ and $\varepsilon_n$ be given by
	Lemma \ref{lem:scale-choice}, and set
	$W_n=U_{R_n,\varepsilon_n}$.
	By \eqref{eq:scale-choice-errors}, we have
	\[
	M(W_n)=\varepsilon_n^2M(U_{R_n})\to0.
	\]
	For $n$ sufficiently large, define
	\[
	H_n=\sqrt{\rho^2-M(W_n)}\,h_{\sigma_n}.
	\]
	Choose $\sigma_n\to-\infty$ such that
	\[
	\operatorname{supp}W_n\cap\operatorname{supp}H_n
	=\varnothing,
	\qquad
	A(H_n)=o(\varepsilon_n^{\gamma_q}).
	\]
	Set $v_n=W_n+H_n$. Their disjoint supports give $M(v_n)=M(W_n)+M(H_n)=\rho^2$, so $v_n\in S_\rho$, and
	\[
	A(v_n)=A(W_n)+A(H_n),
	\qquad
	B_q(v_n)=B_q(W_n)+B_q(H_n).
	\]
We use the positivity of the kernel to obtain
	\[
	\mathcal D(v_n)
	\geq \mathcal D(W_n)+\mathcal D(H_n).
	\]
	Using \eqref{eq:UR-A-D} and the scaling identities, we deduce for all $s\in\mathbb R$ that
	\begin{equation}\label{eq:test-energy-upper}
		E_\beta(s\star v_n)
		\leq{}
		A(U_{R_n})
		\left(\frac{e^{2s}}2-\frac{e^{2ps}}{2p}\right)
		-\frac{\beta e^{r_qs}}q
		\varepsilon_n^{\gamma_q}B_q(U_{R_n})
		+\frac{e^{2s}}2A(H_n).
	\end{equation}
	
	Since $A(U_{R_n})\to A(U)>0$, there exist constants $C_0,c_0>0$
	independent of $n$ such that
	\[
	A(v_n)\leq C_0,
	\qquad
	\mathcal D(v_n)\geq\mathcal D(W_n)
	=A(U_{R_n})\geq c_0.
	\]
	Consequently,
	\[
	E_\beta(s\star v_n)
	\leq \frac{C_0}{2}e^{2s}
	-\frac{c_0}{2p}e^{2ps}.
	\]
	Since $p>1$, we can choose a constant $S>0$ independent of $n$ such that
	\[
	\sup_{|s|\geq S}E_\beta(s\star v_n)
	\leq \frac12T_{\mu,a,b}.
	\]
	
	For $|s|\leq S$, we use
	\[
	\frac{e^{2s}}2-\frac{e^{2ps}}{2p}
	\leq\frac{p-1}{2p},
	\qquad
	e^{r_qs}\geq e^{-r_qS},
	\]
	together with \eqref{eq:truncated-B-lower}.
	It follows from \eqref{eq:test-energy-upper} that
	\[
	\begin{aligned}
		\max_{|s|\leq S}E_\beta(s\star v_n)
		&\leq \frac{p-1}{2p}A(U_{R_n})
		-\frac{\beta c_q}{q}e^{-r_qS}
		\varepsilon_n^{\gamma_q}
		+\frac{e^{2S}}2A(H_n)\\
		&\leq T_{\mu,a,b}+O(R_n^{-\kappa})
		-c\beta\varepsilon_n^{\gamma_q}
		+o(\varepsilon_n^{\gamma_q}),
	\end{aligned}
	\]
	where $c>0$ is independent of $n$.
	Combining this upper bound with the estimate for $|s|\geq S$
	and using $R_n^{-\kappa}=o(\varepsilon_n^{\gamma_q})$,
	we deduce for $n$ sufficiently large that
	\begin{equation}\label{eq:test-max-strict}
		\max_{s\in\mathbb R}E_\beta(s\star v_n)
		\leq T_{\mu,a,b}
		-c\beta\varepsilon_n^{\gamma_q}
		+o(\varepsilon_n^{\gamma_q})
		<T_{\mu,a,b}.
	\end{equation}
	
	By Lemma \ref{lem:critical-super-fibers},
	$t(v_n)\star v_n\in\mathcal M_{\rho,\beta}$ and
	\[
	m_\beta
	\leq E_\beta(t(v_n)\star v_n)
	=\max_{s\in\mathbb R}E_\beta(s\star v_n)
	<T_{\mu,a,b}.
	\]
	This establishes \eqref{eq:absolute-strict-level}.
\end{proof}

\begin{proposition}[Normalized ground state]
	\label{prop:critical-super-ground}
	Assume $q_c\leq q<2^\sharp$, with $0<\beta<\beta_*$
	when $q=q_c$, and assume
	\eqref{eq:absolute-truncation-condition}.
	Then $m_\beta$ is attained by a positive normalized
	ground state, and
	\begin{equation}\label{eq:m-below-T}
		0<m_\beta<T_{\mu,a,b}.
	\end{equation}
	When $q>q_c$, the conclusion holds for every $\beta>0$.
\end{proposition}

\begin{proof}
	We combine Lemma \ref{lem:coercivity-M} with
	Proposition \ref{prop:absolute-strict-level} to obtain
	\eqref{eq:m-below-T}.
	Let $\{u_n\}\subset\mathcal M_{\rho,\beta}$ be the bounded
	Palais--Smale sequence obtained in Lemma \ref{lem:m-PS}.
	Since $P_\beta(u_n)=0$, Theorem \ref{thm:global-compactness-intro}\textit{(ii)}
	applies. In view of $0<m_\beta<T_{\mu,a,b}$, its weak
	limit $u^0$ satisfies $u^0 \not\equiv 0$. The same theorem implies
	$M(u^0)=\rho^2$ and ensures that $u^0$ is a normalized
	solution.
	
	By Proposition \ref{prop:pohozaev}, $P_\beta(u^0)=0$.
	Thus $u^0\in\mathcal M_{\rho,\beta}$ and
	$E_\beta(u^0)\geq m_\beta$.
	The energy decomposition yields
	\[
	m_\beta
	=E_\beta(u^0)+\sum_{j=1}^{k}J(U^j)
	\geq m_\beta+kT_{\mu,a,b}.
	\]
	Since $T_{\mu,a,b}>0$, we have $k=0$. Consequently,
	\[
	u_n\to u^0 \quad\text{in }X,
	\qquad
	E_\beta(u^0)=m_\beta.
	\]
	
	Define $u_\beta=|u^0|$.
	The invariance of $M$, $E_\beta$, and $P_\beta$
	under the mapping $u\mapsto|u|$ implies
	\[
	u_\beta\in\mathcal M_{\rho,\beta},
	\qquad
	E_\beta(u_\beta)=m_\beta.
	\]
	By the uniqueness of the fiber projection and $t(u_\beta)=0$,
	we have
	\[
	\Psi(u_\beta)=E_\beta(u_\beta)=m_\beta
	=\inf_{v\in S_\rho}\Psi(v).
	\]
	Therefore $d\Psi(u_\beta)=0$ and
	Lemma \ref{lem:natural-constraint} implies that $u_\beta$
	is a constrained critical point of $E_\beta$ on $S_\rho$.
	By applying local regularity and the strong maximum principle, we deduce that $u_\beta>0$ in $\mathbb R^N\setminus\{0\}$. Moreover, using Proposition \ref{prop:pohozaev}, we obtain
	\[
	\lambda_\beta\rho^2
	=-\beta(1-\delta_q)B_q(u_\beta)<0.
	\]
	
	Finally, every normalized solution belongs to
	$\mathcal M_{\rho,\beta}$ by
	Proposition \ref{prop:pohozaev}, and therefore has energy
	at least $m_\beta$. Thus $u_\beta$ is a normalized
	ground state.
\end{proof}

\begin{proof}[Proof of Theorem \ref{thm:main}]
We prove part \textit{(i)} using Proposition
\ref{prop:negative-ground-state} and establish part \textit{(ii)} by
applying Propositions \ref{prop:relative-strict-level} and
\ref{prop:second-solution}. For $q=q_c$, we observe that
the truncation condition \eqref{eq:absolute-truncation-condition}
is equivalent to $\kappa>4d/(N+2)$. We then prove part \textit{(iii)}
using Propositions \ref{prop:absolute-strict-level} and
\ref{prop:critical-super-ground}. Finally, we use the same
propositions to establish part \textit{(iv)} for every $\beta>0$.
\end{proof}

\section*{Declarations}

\noindent\textbf{Conflict of interest.} The authors declare that they have no conflict of interest.

\medskip
\noindent\textbf{Data availability.} Data sharing is not applicable to this article as no datasets were generated or analyzed during the current study.

\end{document}